\documentclass[a4paper,10pt]{article}

\usepackage{stmaryrd}
\usepackage{amsfonts}
\usepackage{bbm}

\usepackage{mathrsfs}
\usepackage{latexsym,amssymb,amsmath,amscd,amsthm,amsxtra}
\usepackage[utf8]{inputenc}
\usepackage[T1]{fontenc}
\usepackage{lmodern}
\usepackage{amssymb}
\usepackage[all]{xy}
\usepackage{nicefrac,mathtools,enumitem}
\usepackage{microtype}

\usepackage{color}
\usepackage{xcolor}
\newtheorem{thm}{Theorem}[section]
\newtheorem{lem}[thm]{Lemma}
\newtheorem{cor}[thm]{Corollary}
\newtheorem{pro}[thm]{Proposition}
\newtheorem{ex}[thm]{Example}
\newtheorem{rmk}[thm]{Remark}
\newtheorem{defi}[thm]{Definition}

\newcommand{\be }{\begin{equation}}
\newcommand{\ee }{\end{equation}}

\newcommand{\pf}{\noindent{\bf Proof.}\ }
\newcommand {\emptycomment}[1]{} %to remove paragraphs

\newcommand{\huaA}{\mathcal{A}}%{{\mathcal{F}}}%{\mathcal{A}}

\newcommand{\huaM}{\mathcal{M}}

\newcommand{\CWM}{C^{\infty}(M)}

\newcommand{\frka}{\mathfrak a}

\newcommand{\frkg}{\mathfrak g}

\newcommand{\frkk}{\mathfrak k}

\newcommand{\frkX}{\mathfrak X}

\def\qed{\hfill ~\vrule height6pt width6pt depth0pt}

\newcommand{\half}{\frac{1}{2}}

\newcommand{\pair}[1]{\left( #1\right)_+}
\newcommand{\pairm}[1]{\left( #1\right)_-}
\newcommand{\pairb}[1]{{\{} #1 {\} }}

\newcommand{\Courant}[1]{\left\llbracket  #1\right\rrbracket }
\newcommand{\dorfman}[1]{\left\lceil  #1\right\rceil }

\newcommand{\br}[1]{   [ \cdot,    \cdot  ]   }

\newcommand{\Id}{\rm{id}}

\newcommand{\g}{\mathfrak g}

\newcommand{\id}{\mathbbm{i}}

\newcommand{\dM}{\mathrm{d}}

\newcommand{\Sym}{\mathsf{Sym}}

\newcommand{\gl}{\mathfrak {gl}}

\newcommand{\Ker}{\mathrm{ker}}

\newcommand{\ad}{\mathrm{ad}}

\newcommand{\pr}{\mathrm{pr}}

\newcommand{\Img}{\mathrm{Im}}

\newcommand{\Vect}{\mathsf{Vect}}
\newcommand{\sgn}{\mathrm{sgn}}
\newcommand{\Ksgn}{\mathrm{Ksgn}}

\newcommand{\V}{\mathbb{V}}

\begin{document}
\title{
Strong homotopy Lie algebras, homotopy Poisson manifolds and Courant algebroids\thanks
 {
Research supported by NSFC (11101179,11471139) and NSF of Jilin Province (20140520054JH). Xiaomeng Xu was partially supported by the SNSF grants P2GEP2-165118 and NCCR SwissMAP.
 }}
\author{{ Honglei Lang$^a$,~~Yunhe Sheng$^b$ ~~ and ~~ Xiaomeng Xu$^c$}\\
{$^a$Department of Mathematics, Peking University, Beijing 100871, China}\\
$^b $Department of Mathematics, Jilin University,
 Changchun 130012,  China\\
$\&$ Kavli Institute for Theoretical Physics China, CAS, Beijing 100190, China\\
 $^c$Section of Mathematics, University of Geneva \\
2-4 Rue de Li$\rm\grave{e}$vre, c.p. 64, 1211-Gen$\rm\grave{e}$ve
4, Switzerland.\\
 { email: hllang@pku.edu.cn,\quad shengyh@jlu.edu.cn,\quad
Xiaomeng.Xu@unige.ch} }

%\author[Honglei Lang]{Honglei Lang}
%\address{Department of Mathematics, Peking University, Beijing 100871, China}
%\email{hllang@pku.edu.cn}
%
%\author[Yunhe Sheng]{Yunhe Sheng}
%\address{Department of Mathematics, Jilin University,
% Changchun 130012,  China,
%$\&$ Kavli Institute for Theoretical Physics China, CAS, Beijing 100190, China}
%\email{shengyh@jlu.edu.cn}
%
%\author[X. Xu]{Xiaomeng Xu}
%\address{Section of Mathematics, University of Geneva \\
%2-4 Rue de Li$\rm\grave{e}$vre, c.p. 64, 1211-Gen$\rm\grave{e}$ve
%4, Switzerland}
%\email{xiaomeng.xu@unige.ch}

\date{}
\maketitle
\begin{abstract}

We  study Maurer-Cartan elements on homotopy Poisson manifolds of degree $n$. They unify many twisted or homotopy structures in Poisson geometry and mathematical physics, such as  twisted Poisson manifolds, quasi-Poisson $\g$-manifolds, and twisted Courant algebroids. %As a byproduct, we justify that the homotopy symplectic manifold of degree $n$ is a homotopy version of the symplectic NQ-manifold and the Maurer-Cartan equation is a homotopy version of the master equation.
%\textcolor{red}{Note that the cotangent bundle of a homotopy Poisson manifold of degree $n$ is a symplectic NQ-manifold of degree $n+1$.}
Using the fact that the dual of an $n$-term $L_\infty$-algebra is a homotopy Poisson manifold of degree $n-1$, we obtain a Courant algebroid  from a $2$-term $L_\infty$-algebra $\g$ via the degree $2$ symplectic NQ-manifold $T^*[2]\g^*[1]$. By integrating the Lie quasi-bialgebroid associated to the Courant algebroid, we obtain a Lie-quasi-Poisson groupoid from a $2$-term $L_\infty$-algebra, which is proposed to be the geometric structure on the dual of a Lie $2$-algebra. These results  lead to a construction of a new 2-term $L_\infty$-algebra from a given one, which could produce many interesting examples.

\end{abstract}
%\tableofcontents

\section{Introduction}

 The concept of an $L_\infty$-algebra (sometimes
called a strongly homotopy (sh) Lie algebra)  was originally
introduced in \cite{stasheff:introductionSHLA,stasheff:shla} as a
model for ``Lie algebras that satisfy the Jacobi identity up to all
higher homotopies''. A Lie 2-algebra is a categorification of  a Lie algebra. It is well-known that the category of Lie 2-algebras is equivalent to the category of  2-term
$L_\infty$-algebras
\cite{baez:2algebras}. The structure of a Lie 2-algebra or a $2$-term $L_\infty$-algebra appears in many areas such as
string theory, higher symplectic geometry
\cite{baez:classicalstring,baez:string}, and Courant algebroids
\cite{Roytenberg1}.

A homotopy Poisson algebra is a graded commutative algebra with an $L_\infty$-structure whose brackets satisfy the Leibniz rule.  It has appeared in Voronov's work  \cite{VoronovHigherP,Voronov1} under the name ``higher Poisson structures'' and Cattaneo and Felder's work \cite{AG} under the name ``$P_\infty$-structures''. See also \cite{Bruce, Mehta} for more recent work involving this structure.  A homotopy Poisson algebra structure of degree $n$ on a graded commutative algebra $\frka$  means that there is a homotopy Poisson algebra structure on the shifted space $\frka[n]$.
%The notion of a homotopy Poisson manifold of degree $n$ was given by Mehta in \cite{Mehta}. In particular, homotopy Poisson manifolds of degree $1$ are used to study reduction of symplectic Q-manifolds. Such structures have appeared in the work \cite{Bruce,VoronovHigherP,Voronov1} under the name ``higher Poisson structures'' and Cattaneo and Felder's work \cite{AG} under the name ``$P_\infty$-structures''.
One interesting example of homotopy Poisson manifolds is given by the dual $\g^*[n-1]$ of an $n$-term $L_\infty$-algebra $\g=\g_0\oplus \g_{-1}\oplus\cdots\oplus\g_{1-n}$,
which turns out to be a homotopy Poisson manifold of degree $n-1$. This generalizes the fact that there is a Lie-Poisson structure on the dual space of a Lie algebra to the case of $n$-term $L_\infty$-algebras.

It is known that there is a one-to-one correspondence between Poisson manifolds and symplectic NQ-manifolds of degree $1$, $$(M,\pi)\leftrightsquigarrow (T^*[1]M,Q=\{\pi,\cdot\}),$$ where the bracket $\{\cdot,\cdot\}$ is the canonical Poisson bracket on the cotangent bundle of the manifold $M$. One can further study the cotangent bundle of a homotopy Poisson manifold $\huaM$ of degree $n$.  The shifted cotangent bundle $T^*[n+1]\huaM$ is a symplectic NQ-manifold of degree $n+1$. Then we study Maurer-Cartan elements on  a degree $n$ homotopy Poisson manifold $\huaM$, which are solutions of the Maurer-Cartan equation associated to the $L_\infty$-structure on the function ring $C^\infty(\huaM)$. Recall that there are many types of twisted or homotopy structures in Poisson geometry and mathematical physics.  See \cite{Anton2,strobl, Schwarzbach2,Schwarzbach3,Roytenberg3,Severa,Severa1} for lists of such examples. We provide a unified framework to describe these structures in which they are viewed as solutions of the Maurer-Cartan equation. For example, given a manifold $M$, the cotangent bundle $T^*[1]M$, equipped with the Schouten bracket $l_2=[\cdot,\cdot]_S$, is a homotopy Poisson manifold of degree $1$. Furthermore, a degree 2 function $\pi\in C^{\infty}(T^*[1]M)$ is a Maurer-Cartan element on $T^*[1]M$ if and only if $\pi$ is a Poisson tensor on $M$; in the twisted case, given a closed $3$-form $H$,  the cotangent bundle $T^*[1]M$ with $l_2=[\cdot,\cdot]_S$ and a ternary bracket $l_3=H(\cdot,\cdot,\cdot)$ on $C^\infty(T^*[1]M)$ is a new homotopy Poisson manifold of degree $1$,  and a degree 2 function $\pi$ satisfies the Maurer-Cartan equation  if and only if $\pi$ is a twisted Poisson tensor.

A homotopy Poisson manifold $\huaM$ of degree $n$ is called symplectic if its binary bracket is non-degenerate (See \cite{VoronovHigherP} for the $\mathbb{Z}_2$-graded case). If $\huaM$ is a symplectic N-manifold of degree $n$, then a function $\alpha$ of degree $n+1$, satisfying the Maurer-Cartan equation $l_2(\alpha,\alpha)=0$, induces a differential $Q=l_2(\alpha,\cdot)$ of degree 1 on $\huaM$.  This shows that a differential graded symplectic manifold or symplectic NQ-manifold is a special homotopy symplectic manifold with a Maurer-Cartan element. A general degree $n$ homotopy symplectic manifold $(\huaM,\{l_i\}_{1\leq i< \infty})$ with  $l_1=0$ is therefore a homotopy version of a symplectic N-manifold. In this case, the Maurer-Cartan equation becomes
  $$\frac{1}{2!}l_2(\alpha,\alpha)-\frac{1}{3!}l_3(\alpha,\alpha,\alpha)+\cdots=0,$$
%+\frac{(-1)^{n+2}}{(n+2)!}l_{n+2}(\alpha,...,\alpha)=0$,
which can be viewed as a homotopy version of the classical master equation. In order to relate this to physical applications,  we prove that homotopy symplectic manifolds of degree $n$ with zero $l_1$ are in one-to-one correspondence with twisted symplectic NQ-manifolds with $\Theta|_{\mathcal{M}}=0$ (see Definition \ref{defi:twisted sym NQ}), to which the authors in \cite{IkedaXu} associate sigma models with boundary via AKSZ formalism.

The notion of a Lie bialgebroid was introduced by Mackenzie and Xu \cite{MackenzieXu}. To study the double of a Lie bialgebroid\footnote{There are different ways to describe the double of a Lie bialgebroid, e.g. Mackenzie gave the description of Drinfeld doubles for Lie bialgebroids using  double Lie algebroids in \cite{Mackenzie}; Roytenberg and Voronov gave the description of Drinfeld doubles for Lie bialgebroids using  graded manifolds in \cite{Roytenberg2} and \cite{Voronov0} respectively.}, Liu-Weinstein-Xu introduced the notion of a Courant algebroid \cite{LWXmani}. Then Roytenberg gave an equivalent definition in \cite{Roytenberg4}. Symplectic NQ-manifolds of degree $2$ are in one-to-one correspondence with Courant algebroids \cite{Roytenberg2}. On the other hand, a Courant algebroid gives rise to a $2$-term $L_\infty$-algebra \cite{Roytenberg1}. The present paper uses these facts to construct a new 2-term $L_\infty$-algebra from a given one as follows. Given a Lie 2-algebra $\g=\g_0\oplus \g_{-1}$, its dual $\g^*[1]$ is a homotopy Poisson manifold of degree $1$, which implies that $T^*[2]\g^*[1]$ is a symplectic NQ-manifold of degree $2$. Thus we obtain a Courant algebroid $E=\g_{-1}^*\times (\g_0\oplus\g_0^*)$, where $\g_{-1}^* $ is the base manifold. Choosing some special sections of $E$ and linear functions on the base manifold $\g_{-1}^*$, we obtain a $2$-term $L_\infty$-algebra $\tilde{\g}=\tilde{\g}_0\oplus \tilde{\g}_{-1}$, where $\tilde{\g}_0=\g_0\oplus(\g_{-1}\otimes\g_0^*)$ and $\tilde{\g}_{-1}=\g_{-1}$.  We thus have
$$
\xymatrix{{}\save[]+<1cm,0cm>*\txt<4pc>{%
\mbox{2-term}~ $L_\infty$-algebra ~$\g$}\restore} \longmapsto \xymatrix{{}\save[]+<1cm,0cm>*\txt<4pc>{%
homotopy Poisson manifold $\g^*[1]$}\restore} \longmapsto \xymatrix{{}\save[]+<1cm,0cm>*\txt<6pc>{%
degree 2 symplectic NQ-manifold ~$T^*[2] \g^*[1]$}\restore} \longmapsto\xymatrix{{}\save[]+<1cm,0cm>*\txt<8pc>{%
Courant algebroid $E=\g_{-1}^*\times (\g_0\oplus\g_0^*)$}\restore} \longmapsto \xymatrix{{}\save[]+<1cm,0cm>*\txt<4pc>{%
\mbox{2-term}~ $L_\infty$-\mbox{algebra}~$\tilde{\g}$.}\restore}
$$
In particular, we obtain the 2-term $L_\infty$-algebra associated to the omni-Lie algebra \cite{shengzhu1,Alan} from a vector space $V$, which is viewed as the abelian 2-term $L_\infty$-algebra $(V\stackrel{\Id}{\longrightarrow}V,l_2=0,l_3=0)$, and obtain the 2-term $L_\infty$-algebra of string type \cite{shengzhu1} from  the 2-term $L_\infty$-algebra $(\mathbb R\stackrel{0}{\longrightarrow}\frkk,l_2=[\cdot,\cdot]_\frkk,l_3=0)$, where  $(\frkk,[\cdot,\cdot]_\frkk)$ is a Lie algebra.
%Starting from the string Lie 2-algebra and the 2-term $L_\infty$-algebra $(\frkk\stackrel{\Id}{\longrightarrow}\frkk,l_2=[\cdot,\cdot]_\frkk,l_3=0)$, we can obtain some new interesting examples of 2-term $L_\infty$-algebras. See Examples \ref{ex:1}-\ref{ex:3} for details.
These interesting examples show that the construction of $\tilde{\g}$ from $\g$ has some important properties and applications, which will be studied in the future, since it is unrelated to the key subject of this paper. Nevertheless, we prove that there is a homomorphism from $\tilde{\g}$ to $\g$ (see Theorem \ref{thm:morphism}).

We also observe that the Courant algebroid $E=\g_{-1}^*\times (\g_0\oplus\g_0^*)$ obtained above is the double of a Lie quasi-bialgebroid $(A,\delta,\phi)$, where $A=\g_{-1}^*\times \g_0^*\longrightarrow \g_{-1}^*$ is an action Lie algebroid. By integration, we get a Lie-quasi-Poisson groupoid $(\Gamma,\Pi,\phi)$, where $\Gamma=\mathfrak{g}^*_{-1}\times{\mathfrak{g}^*_0}\rightrightarrows{\mathfrak{g}^*_{-1}}$ is an action groupoid.
In summary, we have
$$
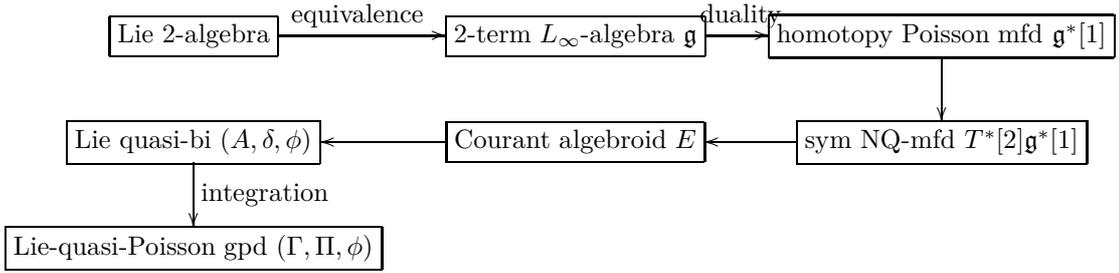
\xymatrix{
*+[F]{\mbox{Lie~ 2-algebra}}\ar[r]^{\mbox{equivalence}\qquad} &*+[F]{\mbox{2-term~} L_\infty\mbox{-algebra}~\g}\ar[r]^{\mbox{duality}~\quad}&*+[F]{\mbox{homotopy Poisson mfd}~\g^*[1]}\ar[d] \\
*+[F]{\mbox{Lie~quasi-bi}~(A,\delta,\phi)}\ar[d]^{\mbox{integration}}&*+[F]{\mbox{Courant~algebroid}~E}\ar[l]&*+[F]{\mbox{sym~NQ-mfd}~T^*[2]\g^*[1]}\ar[l]\\
*+[F]{\mbox{Lie-quasi-Poisson~gpd}~(\Gamma,\Pi,\phi)}&&
 }
$$
If we view the 2-vector space $\mathfrak{g}^*_{-1}\oplus{\mathfrak{g}^*_0}\rightrightarrows{\mathfrak{g}^*_{-1}}$ as the dual of the 2-vector space $\g_0\oplus \g_{-1}\rightrightarrows{\mathfrak{g}_0}$,
this Lie-quasi-Poisson groupoid can be viewed as the natural geometric structure on the ``dual'' of a Lie 2-algebra.

\vspace{2mm}
\noindent {\bf Notations:}
For a graded vector space $V= \sum_{n\in \mathbb Z} V_n $, we use $V[l]$ to denote the $l$-shifted graded vector space, namely
$V[l]_n=V_{l+n}$; we use $\Sym(V)$ to denote the symmetric algebra of $V$. We use $[\cdot,\cdot]_S$ to denote the Schouten bracket of sections of a Lie algebroid.

\section{Preliminaries}\label{Sec:GeneralizedCourant}

Lie algebras can be categorified to Lie 2-algebras.  For a good
introduction on this subject see \cite{baez:2algebras,LadaMarkl}. Vector spaces can be categorified to $2$-vector spaces. Let $\Vect$
be the category of vector spaces.
A {\bf$2$-vector space} is a category in the category $\Vect$.
%One can  make the 2-category of  2-vector spaces into a symmetric monoidal 2-category \cite{street}.
A $2$-vector space $C$ is a category with a vector space of
objects $C_0$ and a vector space of morphisms $C_1$, such that all
the structure maps are linear. Let $s,t:C_1\longrightarrow C_0$ be
the source and target maps respectively.

The 2-category of $2$-vector spaces is
equivalent to the 2-category of 2-term complexes of vector spaces.
Roughly speaking, given a $2$-vector space $C$,
$\Ker(s)\stackrel{t}{\longrightarrow}C_0$ is a 2-term complex.
Conversely, any 2-term complex of vector spaces
$V_1\stackrel{\dM}{\longrightarrow}V_0$ gives rise to a $2$-vector
space of which the set of objects is $V_0$, the set of morphisms is
$V_0\oplus V_1$, the source map $s$ is given by $s(v+m)=v$, and the
target map $t$ is given by $t(v+m)=v+\dM m$, where $v\in V_0,~m\in
V_1.$

A  Lie $2$-algebra is a $2$-vector space $C$ equipped with
 a skew-symmetric bilinear functor,  such that the Jacobi identity is controlled by a natural isomorphism,
which satisfies a coherence law of its own.
The category of  Lie $2$-algebras is  equivalent
to the category of 2-term $L_\infty$-algebras (see \cite{baez:2algebras}).

\begin{defi}
An {\bf  $L_\infty$-algebra} is a graded vector space $\mathfrak{g}=\oplus_{i\in \mathbb Z}\mathfrak{g}_i$ equipped with a system $\{l_k|~1\leq k<\infty\}$
of linear maps $l_k:\wedge^k\mathfrak{g}\longrightarrow \mathfrak{g}$ with degree
$\deg(l_k)=2-k$, where the exterior powers are interpreted in the
graded sense and the following relation with Koszul sign ``{\rm Ksgn}''  is
satisfied for all $n\geq0$:
\begin{equation}\label{eq:higher-jacobi}
\sum_{i+j=n+1}(-1)^{i(j-1)}\sum_{\sigma}\sgn(\sigma)\Ksgn(\sigma)l_j
(l_i(x_{\sigma(1)},\cdots,x_{\sigma(i)}),x_{\sigma(i+1)},\cdots,x_{\sigma(n)})=0.
\end{equation}
Here the summation is taken over all $(i,n-i)$-unshuffles with
$i\geq1$.
\end{defi}

We denote an $L_\infty$-algebra by $(\oplus_{i\in \mathbb Z}\mathfrak{g}_i,\{l_i\}_{i\geq1})$, or simply by $\g$ if there is no confusion. In particular, we denote a 2-term $L_\infty$-algebra by $(\g_{-1}\stackrel{l_1}{\longrightarrow}\g_0,l_2,l_3)$.

\begin{defi}\label{defi:Lie-2hom}
Let $\g=(\g_{-1}\stackrel{l_1}{\longrightarrow}\g_0,l_2,l_3)$ and $\g'=(\g_{-1}'\stackrel{l_1'}{\longrightarrow}\g_0',l_2',l_3')$ be two $2$-term  $L_\infty$-algebras. An {\bf
$L_\infty$-algebra homomorphism} $F$ from $\frkg$ to $ \frkg'$ consists of:
 linear maps $F_0:\frkg_0\rightarrow \frkg_0',~F_1:\frkg_{-1}\rightarrow \frkg_{-1}'$
 and $F_{2}: \frkg_{0}\wedge \frkg_0\rightarrow \frkg_{-1}'$,
such that the following equalities hold for all $ x,y,z\in \frkg_{0},
m\in \frkg_{-1},$
\begin{itemize}
\item [$\rm(i)$] $F_0\circ l_1=l_1'\circ F_1$,
\item[$\rm(ii)$] $F_{0}l_2(x,y)-l_2'(F_{0}(x),F_{0}(y))=l_1'F_{2}(x,y),$
\item[$\rm(iii)$] $F_{1}l_2(x,m)-l_2'(F_{0}(x),F_{1}(m))=F_{2}(x,l_1( m))$,
\item[$\rm(iv)$]
$F_2(l_2(x,y),z)+c.p.+F_1(l_3(x,y,z))=l_2'(F_0(x),F_2(y,z))+c.p.+l_3'(F_0(x),F_0(y),F_0(z))$.
\end{itemize}
\end{defi}

Given  a graded vector space $V= \sum_{n\in \mathbb Z} V_n
$, it is well-known that there is a graded version of Kosmann-Schwarzbach's big bracket \cite{Schwarzbach1}, which we denote by
$\pairb{\cdot,\cdot}$, on $\Sym(V^*[l]) \otimes \Sym(V[k]) \cong
\Sym(V^*[l]\oplus V[k]) \cong \Sym(T^*[l+k]V[k])$ by extending the
usual pairing between $V^*$ and $V$ via a graded Leibniz
rule
\begin{eqnarray} \label{eq:g-leibniz}
\pairb{u, v\wedge w}&=& \pairb{u, v}\wedge w + (-1)^{(|u|+l+k) |v|}
v\wedge \pairb{u, w},\\
\label{eq:g-symmetry} \pairb {u, v}&=&-(-1)^{(|u|+k+l)(|v|+k+l)}
\pairb{v, u},
\end{eqnarray} where $u \in \Sym(V^*[l]\oplus V[k]) _{|u|}$ and $v \in
\Sym(V^*[l]\oplus V[k]) _{|v|}$. The big bracket is in fact the
canonical graded Poisson bracket on $T^*[-l-k]V^*[-k]$. Thus, we have a
graded Jacobi identity:
\begin{equation}\label{eq:g-jacobi}
\pairb{u, \pairb{v, w}}=\pairb{\pairb{u, v},
w}+(-1)^{(|u|+k+l)(|v|+k+l)} \pairb{v,
  \pairb{u, w}}.
\end{equation}

We view $l_i:\wedge^i\mathfrak{g}\longrightarrow \mathfrak{g}$ as elements in $\Sym^i(\g^*[-1])\otimes \g[k]$. Then we have

\begin{lem} \label{lem:delta-big-bracket}
A series of degree $2-k$ elements $l_i\in (\Sym^i(\g^*[-1]) \otimes
\g[k])$ with $i=1, 2, \dots$ on $\g=\g_0\oplus \g_{-1}
\oplus \dots$ gives an
$L_\infty$-algebra structure
 if and only if $\pairb{\sum_{i=1}^\infty l_i,
\sum_{i=1}^\infty l_i}=0$.
\end{lem}

\begin{defi}
A  {\bf Courant algebroid}  is a vector bundle $E\longrightarrow M$, together
with a
fiber metric $\pair{\cdot,\cdot}$ (so we can identify $E$ with
$E^*$), a bundle map $\rho: E \longrightarrow TM$ (called the
anchor), a bilinear bracket operation (Dorfman bracket)
$\dorfman{\cdot,\cdot}$ on $\Gamma(E)$,
such that for all $e_1,e_2,e_3\in\Gamma(E)$, we have
\begin{eqnarray}
\label{nonskewsymmetric}\dorfman{e,e}&=&\half D \pair{e,e},\\
  \label{invariantinner} \rho(e_1)\pair{e_2,e_3}&=&\pair{\dorfman{e_1,e_2},e_3}+\pair{e_2,\dorfman{e_1, e_3}},\\
\label{TCcondition3} \dorfman{e_1,\dorfman{e_2,e_3}}&=&\dorfman{\dorfman{e_1,e_2},e_3}+\dorfman{e_2,\dorfman{e_1,e_3}},
\end{eqnarray}
where the operator $D:\CWM\longrightarrow\Gamma(E)$ is the map defined by $$\pair{e,Df}=\rho(e)f.$$
\end{defi}

One can also use the skew-symmetric Courant bracket $\Courant{\cdot,\cdot}:$
$
\Courant{e_1,e_2}=\half(\dorfman{e_1,e_2}-\dorfman{e_2, e_1}).
$
But the equality \eqref{TCcondition3} does not hold anymore. We have
$$
\Courant{\Courant{e_1,e_2},e_3}+\Courant{\Courant{e_2,e_3},e_1}+\Courant{\Courant{e_3,e_1},e_2}=DT(e_1,e_2,e_3),
$$
where $T(e_1,e_2,e_3)\in\CWM$ is defined by
$$
T(e_1,e_2,e_3)=\frac{1}{6}\pair{\Courant{e_1,e_2},e_3}+c.p..
$$

\begin{thm}\label{thm:CALie2}{\rm(\cite{Roytenberg1})}
  A Courant algebroid structure on a vector bundle $E$ gives rise to a $2$-term $L_\infty$-algebra structure on $\CWM\oplus\Gamma(E)$, where $l_i$ are given by
  \begin{equation}\label{eqn:defil2}
  \left\{\begin{array}{rcll}
  l_1(f)&=&Df,&~\forall ~f\in\CWM,\\
   l_2(e_1,e_2)&=& \Courant{e_1,e_2} & ~\forall~e_1,e_2\in\Gamma(E),\\
  l_2(e_1,f)&=& \half\pair{e_1,Df} &
  ~\forall~e_1\in\Gamma(E),f\in\CWM,\\
  l_3(e_1,e_2,e_3)&=&-T(e_1,e_2,e_3) &~\forall~e_1,e_2,e_3\in\Gamma(E).
   \end{array}\right.
\end{equation}
\end{thm}

A Courant algebroid can be described by a symplectic NQ-manifold of degree $2$ \cite{Roytenberg2}.
Explicitly, let $(E,\pair{\cdot,\cdot})$ be a pseudo-Euclidean vector bundle over a manifold $M$. Then, $E[1]$ is a Poisson N-manifold of degree $-2$. Let $\huaM$ be its minimal symplectic realization. An NQ-structure on $\huaM$ is determined by a cubic Hamiltonian function $\Theta$ on $\huaM$ satisfying $\{\Theta,\Theta\}=0$, where $\{\cdot,\cdot\}$ is the even Poisson bracket on $\huaM$ of degree $-2$. Such a $\Theta$ corresponds to a Courant algebroid structure on $E$. In fact, we define the anchor $\rho$ and the Dorfman bracket $\dorfman{\cdot,\cdot}$ as the derived bracket by
\[\rho(e)f=\{\{e,\Theta\},f\},\quad \quad \dorfman{e_1,e_2}=\{\{e_1,\Theta\},e_2\},\quad  \quad ~\forall~f\in C^\infty(M), ~\forall~e,e_1,e_2\in \Gamma(E).\]
Here  $C^\infty(M)$ and $\Gamma(E)$ are the algebras of  degree $0$ and degree $1$ polynomial functions on $\huaM$ respectively. We refer to \cite{Voronov0} for  a general construction of doubles for graded QS-manifolds and QP-manifolds.

For a vector bundle $A$, consider the graded manifold $T^*[2]A[1]$. It is canonically equipped with a Poisson bracket of degree $-2$, and is actually the minimal symplectic realization of $A\oplus A^*$. This Poisson bracket, called the {\bf big bracket} in \cite{Schwarzbach1}, is denoted here by $\{\cdot,\cdot\}$. Let $(x^i,
\xi^a)$ be local coordinates on $A[1]$, we denote by $(x^i,
\xi^a, \theta_a,p_i)$ the local coordinates on  $T^*[2]A[1]$. About their degrees, we have
$$
\mbox{degree}(x^i,
\xi^a, \theta_a,p_i)=(0,1,1,2).
$$
The big bracket satisfies
$$
\{x^i,p_j\}=\delta^i_j=-\{p_j,x^i\},\quad \{\xi^a,\theta_b\}=\delta^i_j=\{\theta_b,\xi^a\}.
$$
A Lie algebroid structure on $A$ is equivalent to a degree $3$ function $\mu=\rho^i_bp_i\xi^b+\half\mu^a_{bc}\xi^b\xi^c\theta_a$ such that $\{\mu,\mu\}=0.$

A {\bf Lie quasi-bialgebroid} structure on $A$ is given by a degree $3$ function $\mu+\gamma+\phi$, which can be locally written as
$$
\mu=\rho^i_bp_i\xi^b+\half\mu^a_{bc}\xi^b\xi^c\theta_a,\quad \gamma=\varrho^{ib}p_i\theta_b+\half\gamma_a^{bc}\xi^a\theta_b\theta_c,\quad \phi=\frac{1}{6}\phi^{abc}\theta_a\theta_b\theta_c,
$$
and they satisfy
$$
\{\mu+\gamma+\phi,\mu+\gamma+\phi\}=0.
$$
More precisely, we have
\begin{equation}\label{eq:defiquasi}
\{\mu,\mu\}=0,\quad \{\mu,\gamma\}=0,\quad \half\{\gamma,\gamma\}+ \{\mu,\phi\}=0,\quad\{\gamma,\phi\}=0.
\end{equation}
See  \cite{Schwarzbach2} for more details.
On $A\oplus A^*$, there is a natural Courant algebroid structure, in which the degree $3$ function $\Theta$ is exactly $\mu+\gamma+\phi$. Note that there are two canonical fiber metric $(\cdot,\cdot)_{\pm}$ on $A\oplus A^*$:
 \begin{equation}
  (x+\xi,y+\eta)_{\pm}=\langle\xi,y\rangle\pm\langle x,\eta\rangle,
 \end{equation}
where $\langle\cdot,\cdot\rangle$ is the usual pairing without any degree involved in.

\section{Maurer-Cartan elements on homotopy Poisson manifolds}
%The notions of a homotopy Poisson algebra (manifold) of degree $n$ was given in \cite{Mehta}. Here we recall them briefly.
First we recall the notions of homotopy Poisson algebras and homotopy Poisson manifolds. They are also called  $P_\infty$-manifolds in \cite{AG,Schatz} and higher Poisson manifolds in \cite{Bruce,VoronovHigherP,Voronov1}. Here we follow the convention in \cite{Mehta}.

\begin{defi}
A {\bf  homotopy Poisson algebra} of degree $n$   is a graded commutative algebra $\frka$ over a field of characteristic zero with an $L_\infty$-algebra structure $\{l_m\}_{m\geq1}$ on $\frka[n]$, such that the map
\[x\longrightarrow{l_m(x_1,\cdots,x_{m-1},x)},\ \ \ \ x_1,\cdots,x_{m-1},x\in \frka\]
is a derivation of $\frka$ of degree $\kappa:=2-m-n(m-1)+\sum_{i=1}^{m-1}|x_i|$, i.e., for all $x,y\in \frka$,
\begin{eqnarray*}
l_m(x_1,\cdots,x_{m-1},xy)=l_m(x_1,\cdots,x_{m-1},x)y+(-1)^{\kappa|x|}
x l_m(x_1,\cdots,x_{m-1},y).
\end{eqnarray*}
Here, $|x|$ denotes the degree of $x\in\frka.$
%\comment{here I change '$+\sum_{i=1}^{m-1}|x_i|$' to $-\sum_{i=1}^{m-1}|x_i|$ in order to be coherent with Rmk 3.2 (i), is  it Ok?}\comment{I don’t think it is good. Since the degree of $x_i$ here is $|x_i|-n$, take sum of them and plus the original $2-m$ we get $2-m-n(m-1)+\sum_{i=1}^{m-1}|x_i|$. Although we can understand $1-(n+1)(m-1)+\sum_{i=1}^{m-1}|x_i|$ from $a[1]$. But in our definition of $L_\infty$-algebras, we only use $a$. So I suggest we use $2-m-n(m-1)+\sum_{i=1}^{m-1}|x_i|$.}

A  homotopy Poisson algebra of degree $n$ is of {\bf finite type} if there exists a $q$ such that $l_m=0$ for all $m>q$.

A {\bf homotopy Poisson manifold} of degree $n$ is a graded manifold $\huaM$ whose algebra of functions $C^\infty(\huaM)$ is equipped with a degree $n$ homotopy Poisson algebra structure of finite type.
\end{defi}

Throughout this paper, we use $(\mathcal{M},\{l_i\}_{1\leq i< \infty})$ to denote a homotopy Poisson manifold. Obviously, a usual Poisson manifold is a homotopy Poisson manifold of degree $0$.

\begin{rmk}
  %In this remark, we give several related structures.
  In this remark, we compare the related structures in the literature.

  \begin{itemize}
    \item[\rm(i)] In \cite{AG}, the authors introduced the notion of
a {\bf  $P_\infty$-algebra},  which is a graded commutative algebra $\frka$  with an $L_\infty$-algebra structure  such that the Leibniz rule is satisfied. So a $P_\infty$-algebra is a homotopy Poisson algebra of degree $0$.

\item[\rm(ii)] The notion of {\bf a higher Poisson structure} was introduced in \cite{Voronov1}, and further studied in \cite{Bruce,VoronovHigherP}, where the authors used the superized version of an $L_\infty$-algebra, i.e., $\mathbb{Z}_2$-graded.

\item[\rm(iii)] A {\bf  graded Poisson algebra of degree $k$} is a graded commutative algebra $\frka$ with a degree $-k$ Lie bracket, such that the bracket is a biderivation  of the product, namely
\[[x,y\cdot z]=[x,y]\cdot z+(-1)^{|y|(|x|+k)}y\cdot [x,z].\]
See \cite{ADR} for more details.
Thus, a graded Poisson algebra of degree $k$ is  a homotopy Poisson algebra of degree $k$. In particular, the associated $L_{\infty}$-algebra has only one non-zero map $l_2$.
  \end{itemize}
\end{rmk}

In the following, we list a few interesting examples of different kinds.

\begin{ex}\rm{
Let $A$ be a Lie algebroid. Consider its dual vector bundle, $A^*[1]$, which is  an N-manifold of degree $1$. Its algebra of polynomial functions is
$$
\cdots\oplus \Gamma(A)\oplus C^\infty(M),
$$
where $\Gamma(A)$ is of degree 1, and $C^\infty(M)$ is of degree 0. The Poisson bracket  is in fact the Schouten bracket $[\cdot,\cdot]_S$ on $\Gamma(\wedge^\bullet A)$. It is straightforward to see that $A^*[1]$ is a homotopy Poisson manifold of degree $1$.}
\end{ex}

\begin{ex}\label{ex:cotangent bundle}\rm{For an arbitrary manifold $M$, the shifted cotangent bundle $T^*[1]M$ is a symplectic N-manifold of degree $1$ (see \cite{CZ,Roytenberg2}).
%The symplectic structure is the canonical one on the cotangent bundle.
Its algebra of polynomial functions  is
 $$\cdots\oplus \frkX(M)\oplus C^\infty(M),$$
  where $\frkX(M)$ is of degree 1, and $C^\infty(M)$ is of degree 0. The Poisson bracket is exactly the Schouten bracket $[\cdot,\cdot]_S$ on $\Gamma(\wedge^\bullet TM)$. It is straightforward to see that $(T^{*}[1]M,l_2=[\cdot,\cdot]_S)$ is a homotopy Poisson manifold of degree $1$.
 Similarly, any symplectic N-manifold of degree $n$ is a homotopy Poisson manifold of degree $n$. The $L_\infty$-algebra structure here only contains one nonzero $l_2$.}
\end{ex}

\begin{ex}\label{ex:dual Lie n}\rm{
Given a $2$-term $L_\infty$-algebra $\g=(\mathfrak{g}_ {-1}\stackrel{l_1}\longrightarrow{\mathfrak{g}_0},l_2,l_3)$, its graded dual space $\g^*[1]=\mathfrak{g}^*_0[1]\oplus{\mathfrak{g}^*_ {-1}[1]}$ is an N-manifold of degree $1$ with the base manifold $\g_{-1}^*$.  Its algebra of polynomial functions is
$$
\cdots\oplus \big(C^\infty(\g_{-1}^*)\otimes\g_{0}\big)\oplus C^\infty(\g_{-1}^*).
$$
There is a degree $1$ homotopy Poisson algebra structure on it obtained by extending the original $2$-term $L_\infty$-algebra structure using the Leibniz rule. Thus, the dual of a $2$-term $L_\infty$-algebra is a homotopy Poisson manifold of degree $1$. This generalizes the fact that the dual of a Lie algebra is a linear Poisson manifold.

Similarly, the dual of an $n$-term $L_\infty$-algebra is a homotopy Poisson manifold of degree $n-1$. }
\end{ex}

\begin{ex}{\rm
Given a splitting Lie $n$-algebroid\footnote{ By definition, there is an $L_\infty$-algebra structure on the complex of section spaces together with some compatibility conditions. See \cite{Poncin,shengzhu} for more details. }%In \cite{Bruce}, the author use the terminology of an $L_\infty$-algebroid and study its relation with higher Schouten/Poisson structures. However, all the objects therein are $\mathbb Z_2$-graded. Thus, it is not the same as our Lie $n$-algebroids.
 $\huaA=A_0\oplus{\cdots}\oplus A_{1-n}$, its dual $\huaA^*[n]=A_0^*[n]\oplus{\cdots}\oplus A_{1-n}^*[n]$ is an N-manifold of degree $n$. Its algebra of polynomial functions is
\[\cdots\oplus\big(\Gamma(A_0)\oplus{\Gamma(\wedge^n{A_{1-n}})}\oplus\cdots\big)\oplus\cdots\oplus{\Gamma(A_{1-n})}\oplus{C^{\infty}(M)}.\]
%\comment{degree n can be the composition of many cases, not only these two terms, like $\wedge^{n-1} A_{1-n}\wedge A_{2-n}$}
 There is a degree $n$ homotopy Poisson algebra structure on it. Thus, the dual of a splitting Lie $n$-algebroid is a homotopy Poisson manifold of degree $n$.
}
\end{ex}

It is known that $(M,\pi)$ is a Poisson manifold if and only if $(T^*[1]M,Q)$ is a symplectic NQ-manifold of degree $1$, where the homological vector field $Q$ is given by $Q=\{\pi,\cdot\}$. More generally, the cotangent bundle of a homotopy Poisson manifold of degree $n$ gives rise to a symplectic NQ-manifold of degree $n+1$. Consider the homotopy Poisson manifold given in Example \ref{ex:dual Lie n}, we have

\begin{pro}\label{thm:dual Lie n}
Given an $n$-term $L_\infty$-algebra $\mathfrak{g}=({\mathfrak{g}_{0}}\oplus{\cdots}\oplus\mathfrak{g}_{-n+1},\{l_i\}_{1\leq i\leq n+1})$, the cotangent bundle $T^*[n]\mathfrak{g}^{*}[n-1]$ is a symplectic NQ-manifold of degree $n$, where the degree $1$ homological vector field $Q$ is given by
\begin{equation}\label{eq:Q}
Q=\{\sum l_i,\cdot\},
\end{equation}
in which $\{\cdot,\cdot\}$ is the canonical Poisson structure, and $\sum l_i\in \Sym(\g^*[-1])\otimes\g[1-n]$ is viewed as a polynomial function of degree $n+1$ on $T^*[n]\mathfrak{g}^{*}[n-1]$.
\end{pro}

Similar to \cite{Marco}, a {\bf Maurer-Cartan element} on a degree $n$ homotopy Poisson manifold $\huaM$ is given by a function $\alpha$ on $\mathcal{M}$ satisfying the Maurer-Cartan equation
\begin{equation}
  \sum_i\frac{(-1)^i}{i!}l_i(\alpha,\cdots,\alpha)=0.
\end{equation}
%we introduce the notion of a Maurer-Cartan element on a homotopy Poisson manifold of degree $n$.
%\begin{defi}
%A {\bf Maurer-Cartan element $\alpha$} on a degree $n$ homotopy Poisson manifold $\mathcal{M}$ is a function on $\mathcal{M}$ satisfying the Maurer-Cartan equation
%\begin{equation}
%  \sum_i\frac{(-1)^i}{i!}l_i(\alpha,\cdots,\alpha)=0.
%\end{equation}
%\end{defi}
\begin{ex}{\rm ({\bf quasi-Poisson $\g$-manifolds})\quad
Let $M$ be a manifold and $(\frkk,[\cdot,\cdot]_\frkk,K)$ a quadratic Lie algebra. Let $\{e_a\}$ be an orthogonal basis of $\frkk$ with respect to the metric $K$, i.e. $K(e_a,e_b)=\delta_{ab}$, and $\{e^a\}$  its dual basis. Then $K$ induces an isomorphism, which we denote by $K^\sharp$, from $\frkk^*$ to $\frkk$ via $K(K^\sharp(\xi),u)=\langle \xi,u\rangle$. More precisely, we have $K^\sharp(e^a)= e_a.$ We define $R\in\wedge^3 \frkk^*$ by $$R(u,v,w)= K([u,v]_\frkk,w), \quad\forall ~u,v,w\in \frkk.$$ Then $\mathcal{M}:=T^*[1]M\times \frkk[1]$ is a homotopy Poisson manifold of degree $1$. In fact, we can define an $L_{\infty}$-algebra structure on $C^{\infty}(\mathcal{M})$ generated by
$$
l_2(X,f)=Xf, \ \ l_2(X,Y)=[X,Y]_S, \ \ l_1(\xi)=\delta(\xi),\ \ l_3(\xi,\eta,\gamma)=K^\sharp(R)(\xi,\eta,\gamma),$$
where $f\in \CWM, X,Y\in\frkX(M)$, $\xi,\eta,\gamma\in \frkk^*,$ and $\delta:\wedge^\bullet \frkk^*\longrightarrow \wedge^{\bullet+1} \frkk^*$ is the coboundary operator associated to  the Lie algebra  $\frkk$.

%Choose a local coordinate $(x_i,p_i,e_a)$ on $\mathcal{M}$, where $\{e_a\}$ is an orthogonal basis of $\g$
%such that $K(e_a,e_b)=\delta^a_b$. Then $K(R)=\frac{1}{6}\sum_{a,b,c}C^c_{ab}e_a\wedge e_b\wedge e_c\in\wedge^3\g$,
%\comment{to make sure $R(e_a,e_b,e_c)=C^c_{ab}$, it needs to add $\frac{1}{6}$ to the right hand of the equality?}
%where $\{C^c_{ab}\}$ are the structure constants of $\g$.
A degree 2 function  $\alpha=\pi+\rho$, where $\pi\in\wedge^2\frkX(M)$ and $\rho\in \frkk^*\otimes \frkX(M)$, is a Maurer-Cartan element, i.e.
$-l_1(\alpha)+\frac{1}{2}l_2(\alpha,\alpha)-\frac{1}{3!}l_3(\alpha,\alpha,\alpha)=0$, if and only if the following three conditions hold:
$$l_1(\rho)=\frac{1}{2}[\rho,\rho]_S,\quad [\pi,\rho]_S=0,\quad \frac{1}{2}[\pi,\pi]_S=\frac{1}{6}K^\sharp(R)(\rho,\rho,\rho).$$
These conditions are equivalent to that $\rho:\frkk\longrightarrow \frkX(M)$ is a Lie algebra morphism, $\pi$ is $\frkk$-invariant and $\frac{1}{2}[\pi,\pi]_S=\wedge^3\rho(K^\sharp(R))$ respectively.
Therefore, a quasi-Poisson $\g$-manifold \cite{Anton1, Anton2} gives rise to a Maurer-Cartan element on $\huaM$.}
\end{ex}

A homotopy Poisson manifold $(\mathcal{M},\{l_i\}_{1 \leq i<\infty})$ is called {\bf symplectic} if the binary bracket $l_2$ is non-degenerate.
%given by a nondegenerate $2$-form $\omega$ on $\mathcal{M}$, i.e. $l_2(f,g)=\omega^{-1}(df,dg)$, for all $f,g\in C^{\infty}(\mathcal{M})$.
We refer the reader to  \cite{VoronovHigherP} for a thorough discussion of homotopy symplectic structures  in the setting of $\mathbb{Z}_2$-graded manifolds.

We turn to study  Maurer-Cartan elements on a homotopy symplectic manifold of degree $n$. First we present some examples with mathematical and physical interests.

\begin{ex}\label{ex:twisted Poisson}{\rm ({\bf twisted Poisson structures})
The shifted cotangent bundle $T^*[1]M$ of a manifold $M$ is canonically a symplectic $N$-manifold of degree $1$. Slightly different from Example \ref{ex:cotangent bundle}, with a choice of a closed $3$-form $H$, we introduce a nontrivial $l_3$ on the algebra of functions of $T^*[1]M$ by $l_3(X,Y,Z)=H(X,Y,Z),$ for all $X,Y,Z\in \frkX(M)$.
The compatibility of $l_2$ and $l_3$ is due to the fact that $H$ is closed. Thus, $(T^*[1]M,l_2=[\cdot,\cdot]_S,l_3=H)$ is a homotopy symplectic manifold of degree $1$. The idea of adding a closed $3$-form to obtain a new $L_\infty$-algebra was first introduced in \cite{Roytenberg5}.

%Choose a local coordinate $(x^i,p_i)$ on $T^*[1]M$.
A degree $2$ function %$\pi=\frac{1}{2}\pi^{ij}(x)p_ip_j$
$\pi$ is a Maurer-Cartan element of $T^*[1]M$ if and only if
\[\frac{1}{2}l_2(\pi,\pi)-\frac{1}{3!}l_3(\pi,\pi,\pi)=0,\]
which is equivalent to $\frac{1}{2}[\pi,\pi]=\wedge^3 {\pi}^\sharp H$, that is, $\pi$ is a twisted Poisson structure (\cite{Roytenberg3,Severa1}) on $M$. }
\end{ex}

\begin{ex}{\rm({\bf twisted Courant algebroids})
  Let  $E\longrightarrow M$ be a vector bundle with a fiber metric $\pair{\cdot,\cdot}$, and $H$ a closed $4$-form on $M$. Let $\huaM$ be its minimal symplectic realization. See \cite{Roytenberg2} for details. Then $\huaA^2$ is the section space of a vector bundle $\mathbb A$, which fits in the following exact sequence:
   $$
   0\longrightarrow \wedge^2E^*\longrightarrow \mathbb A\stackrel{a}{\longrightarrow} TM\longrightarrow 0.
   $$Similar to the treatment in Example \ref{ex:twisted Poisson}, we can  define a new degree $2$ homotopy   Poisson algebra structure on the algebra of functions of $\huaM$ by adding
   $$
   l_4(\chi_1,\chi_2,\chi_3,\chi_4)=H\big(a(\chi_1),a(\chi_2),a(\chi_3),a(\chi_4)\big),\quad \forall~\chi_i\in \huaA^2.
   $$

    Choose a local coordinate $(x^i,p_i,\xi^a)$.  A degree $3$ function $\alpha=\rho^i_ap_i\xi^a-\frac{1}{3!}f_{abc}\xi^a\xi^b\xi^c$ is a Maurer-Cartan element if and only if $\frac{1}{2}l_2(\alpha,\alpha)+\frac{1}{24}l_4(\alpha,\alpha,\alpha,\alpha)=0.$
 Define an anchor $\rho:E\rightarrow TM$ and a derived bracket $\dorfman{\cdot,\cdot}$ on $\Gamma(E)$ by
$$
\rho(e)f=l_2(l_2(e,\alpha),f), \ \
\dorfman{e_1,e_2}=l_2(l_2(e_1,\alpha),e_2).
$$
According to  \cite[Theorem 4.5]{Roytenberg2}, we have
\begin{eqnarray*}
l_2(\alpha,\alpha)&=&\pair{\rho^*(dx^i),\rho^*(dx^j)}p_ip_j+\xi^a\xi^bdx^j\big([\rho(e_a),\rho(e_b)]_S-\rho(\dorfman{e_a, e_b})\big)p_j
\\ &&+\frac{1}{12}\pair{\dorfman{\dorfman{e_a,e_b}, e_c}+\dorfman{e_b,\dorfman{e_a,e_c}}-\dorfman{e_a,\dorfman{e_b,e_c}},e_d}\xi^a\xi^b\xi^c\xi^d.
\end{eqnarray*}
On the other hand, a straightforward calculation gives that $$l_4(\alpha,\alpha,\alpha,\alpha)=\rho^*H.$$
Thus,  the condition for that $\alpha$ is a Maurer-Cartan element is equivalent to   that $(E,K,\rho,\dorfman{\cdot,\cdot},H)$
is a twisted Courant algebroid, which arises from the study of three dimensional sigma models with Wess-Zumino
term \cite{strobl}. See \cite{XuM} for its close relation with coisotropic Cartan geometry. }
\end{ex}

We have seen that the study of a general homotopy symplectic manifold of degree $n$ and its Maurer-Cartan elements is inspired by many interesting geometric structures. Next, we will see that these structures naturally appear in the study of topological field theory.
To be precise, it is known that there exists a systematic method to construct a topological sigma model from a symplectic NQ-manifold, which is called the Alexandrov-Kontsevich-Schwartz-Zaboronsky (AKSZ) formalism \cite{AKSZ}. We will show that this construction can be extended to the case of a degree $n$ homotopy symplectic manifold $(\huaM,\{l_i\}_{2\leq i<\infty})$ with a degree $n+1$ Maurer-Cartan element $\alpha$.
To do this, we need the notion of a twisted symplectic NQ-manifold, which was introduced in \cite{IkedaXu} to describe the structure on the target supermanifold of the AKSZ sigma model with boundary.

\begin{defi}\label{defi:twisted sym NQ}{\rm \cite{IkedaXu}}
Let $(\mathcal{M},\{\cdot,\cdot\}_s)$ be a symplectic $N$-manifold of degree n and $\alpha$ a degree $n+1$ function on it. Then $(\mathcal{M},\{\cdot,\cdot\}_s,\alpha)$ is called a twisted symplectic $NQ$-manifold if there is a degree $n+1$ symplectic $NQ$-manifold $(T^*[n+1]\mathcal{M},\{\cdot,\cdot\},\Theta)$ such that
\begin{itemize}
\item[\rm(a)] $\{\cdot,\cdot\}_s=\{\cdot,\{\cdot,\Theta\}\}|_\mathcal{M}$;

\item[\rm (b)] it satisfies the canonical transformation equation $e^{-\alpha}\Theta|_\mathcal{M}=0$,
where
$$
%e^{-\alpha}\Theta=\Theta-\{\Theta,\alpha\}+\frac{1}{2}\{\{\Theta,\alpha\},\alpha\}-\cdots.
e^{-\alpha}\Theta=\Theta-\{\alpha,\Theta\}+\frac{1}{2}\{\alpha,\{\alpha,\Theta\}\}-\cdots.
$$
\end{itemize}
\end{defi}
The following theorem states that a homotopy symplectic manifold with a Maurer-Cartan element naturally gives rise to a twisted symplectic NQ-manifold.
\begin{thm}
Degree $n$ homotopy symplectic manifolds  $(\huaM,\{l_i\}_{2\leq i<\infty})$ with a degree $n+1$ Maurer-Cartan element $\alpha$ are in one-to-one correspondence
with twisted symplectic NQ-manifolds $(\huaM,\{\cdot,\cdot\}_s,\alpha)$ with $\Theta|_{\huaM}=0$.
\end{thm}
\pf Let  $(\huaM,\{l_i\}_{2\leq i<\infty})$ be a homotopy symplectic manifold of degree $n$  and $\alpha$ a degree $n+1$ Maurer-Cartan element. Then $(T^*[n+1]\mathcal{M},\{\cdot,\cdot\},\Theta=\sum l_i)$  is a symplectic NQ-manifold of  degree $n+1$, and the following equality holds:
 \begin{equation}\label{eq:temp}
l_k(a_1,\cdots,a_k)=\{a_k,\cdots,\{a_2,\{a_1,\sum l_i\}\}\cdots\}|_\mathcal{M}.
\end{equation}
Define $\{\cdot,\cdot\}_s=l_2$, which is a Poisson bracket determined by a nondegenerate closed $2$-form. Now we extend $\alpha$ to a function on $T^*[n+1]\mathcal{M}$, which is constant along the fiber.
By \eqref{eq:temp}, we have $\Theta|_\mathcal{M}=\{\alpha,\Theta\}|_\mathcal{M}=0$ and
$$
%e^{\alpha} \Theta|_\mathcal{M}=\frac{1}{2}\{\{\Theta,\alpha\},\alpha\}|_\mathcal{M}-\frac{1}{6}
%\{\{\{\Theta,\alpha\},\alpha\},\alpha\}|_\mathcal{M}+\cdots=\frac{1}{2}l_2(\alpha,\alpha)
%-\frac{1}{3!}l_3(\alpha,\alpha,\alpha)+\cdots.
e^{-\alpha} \Theta|_\mathcal{M}=\frac{1}{2}\{\alpha,\{\alpha,\Theta\}\}|_\mathcal{M}-\frac{1}{6}
\{\alpha,\{\alpha,\{\alpha,\Theta\}\}\}|_\mathcal{M}+\cdots=\frac{1}{2}l_2(\alpha,\alpha)
-\frac{1}{3!}l_3(\alpha,\alpha,\alpha)+\cdots.
$$
Thus, the condition for that $\alpha$ satisfies  the Maurer-Cartan equation on $(\mathcal{M},\{l_i\}_{2\leq i<\infty})$ is equivalent to that the canonical transformation equation $e^{-\alpha}\Theta|_{\mathcal{M}}=0$ holds on $T^*[n+1]\mathcal{M}$.  \qed\vspace{3mm}

One immediate consequence of the main construction in \cite{IkedaXu} is that, associated to a degree $n$ homotopy symplectic manifold $(\huaM,\{l_i\}_{2\leq i<\infty})$ with a degree $n+1$ Maurer-Cartan element $\alpha$, there is an AKSZ sigma model with boundary. This can be viewed as a generalization of AKSZ formalism to the setting of homotopy symplectic manifolds. In particular, if $(\huaM,\{l_i\}_{2\leq i<\infty})$ is a symplectic NQ-manifold, i.e. $l_i=0$ for $i\geq 3$, it becomes the AKSZ construction.
%Let $(\Sigma,D,\mu)$ be a differential graded $(dg)$ manifold $\Sigma$ with a D-invariant nondegenerate measure
%$\mu$, where $D$ is a differential on $\Sigma$. Let $\mathcal{M}$ be a Poisson$[n,n]$-manifold with a
%zero $l_1$, and $\alpha$ is a Maurer-Cartan element. Then we know $l_2$ induces a
%symplectic form $\omega$ of degree n on $\mathcal{M}$. A symplectic structure on $Map(\Sigma,\mathcal{M})$ is
%constructed by $\omega:=\mu_*ev^*\omega$, where $ev:\Sigma\times Map(\Sigma,\mathcal{M})\rightarrow \mathcal{M}$
%is an evaluation map defined as$$ev:(z,f)\rightarrow f(z),$$
%and the chain map $\mu_*:\Omega(\Sigma\times Map(\Sigma,\mathcal{M}))\rightarrow \Omega(Map(\Sigma,\mathcal{M}))$ is defined by$$\mu_*\omega(y)(v_1,…,v_k)
%=\int_{\Sigma}\mu(x)w(x,y)(v_1,…,v_k),$$here $v$ is a vector field on $\Sigma$.

%In order to construct the action function, we assume that $X$ is a manifold such that $\partial X=\Sigma$. We
%take a fundamental for $\theta$ for the symplectic structure on $\mathcal{M}$ such that $\omega=-d\theta$.
%$$S=i_{\hat{D_{\Sigma}}}{\mu_{\Sigma}}_*ev^*\theta-{\mu_{\Sigma}}_*ev^*\alpha+$$

\section{$2$-term $L_\infty$-algebras and Courant algebroids}
In this section, we obtain a Courant algebroid from a $2$-term $L_\infty$-algebra. Consequently, we give a general construction of a $2$-term $L_\infty$-algebra from a $2$-term $L_\infty$-algebra, which could produce many interesting examples including the $2$-term $L_\infty$-algebra associated to an omni-Lie algebra and the $2$-term $L_\infty$-algebra of string type.

Let $\mathfrak{g}=(\g_{-1}\stackrel{l_1}{\longrightarrow}\g_0,l_2=l^0_2+l^1_2,l_3)$ be a $2$-term $L_\infty$-algebra, where  $$l_1\in{\mathfrak{g}^*_{-1}\otimes{\mathfrak{g}_0}},\quad l^0_2
\in{\wedge^{2}\mathfrak{g}^*_{0}\otimes{\mathfrak{g}_0}}, \quad l^1_2\in{\mathfrak{g}^*_{0}\wedge{\mathfrak{g}^*_{-1}}\otimes{\mathfrak{g}_{-1}}}, \quad l_3
\in{\wedge^{3}\mathfrak{g}^*_{0}\otimes{\mathfrak{g}_{-1}}}.$$ By Proposition \ref{thm:dual Lie n}, the cotangent bundle $T^*[2]\g^*[1]$ is a symplectic NQ-manifold of degree $2$. Now the canonical Poisson structure $\{\cdot,\cdot\}$ is given by \eqref{eq:g-leibniz}-\eqref{eq:g-jacobi}, in which $k=l=-1$. We give the relation between $l_i$ and $\{\cdot,\cdot\}$ by the following lemma.

\begin{lem}\label{lem:relation}
 For all $x,y\in\g_0$ and $m\in\g_{-1}$, we have
\begin{eqnarray*}
  l_1(m)&=&\{m,l_1\}=-\{l_1,m\},\\
  l_2^0(x,y)&=&\{y,\{x,l_2^0\}\}=-\{\{x,l_2^0\},y\},\\
    l_2^1(x,m)&=&\{m,\{x,l_2^1\}\}=-\{\{x,l_2^1\},m\},\\
    l_3(x,y,\cdot)&=&\{y,\{x,l_3\}\}=-\{\{x,l_3\},y\}.
\end{eqnarray*}
\end{lem}

On the other hand,  symplectic NQ-manifolds of degree $2$ are in one-to-one
correspondence with Courant algebroids \cite{Roytenberg2}. Thus, from $T^*[2]\mathfrak{g}^{*}[1]$, we obtain a Courant algebroid $E$:

\begin{equation}\label{eq:E}
 E=
\g_{-1}^*\times (\g_0^*\oplus \g_0)\longrightarrow
 \g_{-1}^*,
\end{equation}
in which the anchor and the Dorfman bracket are defined by the derived bracket using the degree $3$ function $l=-\sum l_i$. Note that $E$ is the direct sum of two vector bundles $A$ and $A^*$, where $A=\g_{-1}^*\times\g_0^*\longrightarrow\g_{-1}^*$ and $A^*=\g_{-1}^*\times\g_0\longrightarrow\g_{-1}^*$. The fiber metric is given by the canonical pairing between $A$ and $A^*$.

It is necessary to give the precise structures on the Courant algebroid $E$.

 \begin{pro}\label{pro:CA structure}
Consider the Courant algebroid $E$ given above. For constant sections $x,y\in{\mathfrak{g}_0}, \xi,\eta\in{\mathfrak{g}_{0}^*}$ and a linear function $m\in{\mathfrak{g}_{-1}}$, we have
\begin{itemize}
  \item[\rm(i)] the anchor of $x$ is a linear vector field, more precisely, $\rho(x)(m)=l_2^1(x,m);$
    \item[\rm(ii)] the anchor of $\xi$ is a constant vector field, more precisely, $\rho(\xi)=-l_1^*(\xi),$ where $l_1^*$ is defined by $\langle l_1^*(\xi),m\rangle=-\langle\xi,l_1(m)\rangle;$
    \item[\rm(iii)] the image of a linear function under the operator $D$ is not a constant section, we have \begin{equation}\label{eq:D}
        Dm=l_{1}(m)-l^1_2(m,\cdot)\in\mathfrak{g}_{0}+\g_0^*\otimes\mathfrak{g}_{-1};
        \end{equation}
    \item[\rm(iv)] the constant sections in $\g_0$ under the Dorfman bracket are not closed, but we have $$\dorfman{x,y}=l^0_2(x,y)+l_3(x,y,\cdot)\in\mathfrak{g}_{0}+\g_0^*\otimes\mathfrak{g}_{-1};$$
      \item[\rm(v)] the Dorfman bracket of two constant sections in $\g_0^*$ is zero, i.e. $\dorfman{\xi,\eta}=0;$
        \item[\rm(vi)] the Dorfman bracket of a constant section in $\g_0^*$ and a constant section in $\g_0$ is a constant section in $\g_0^*$,
        $$\dorfman{x,\xi}=-\dorfman{\xi,x}=l_2^0(x,\cdot)^*\xi\in{\mathfrak{g}^*_0},$$
        where $l_2^0(x,\cdot)\in\gl(\g_0)$ can be viewed as the adjoint map, and $l_2^0(x,\cdot)^*$ is its dual map.%, i.e. $\langle l_2^0(x,\cdot)^*\xi,y\rangle=-\langle\xi,l_2^0(x,y)\rangle.$
\end{itemize}
 \end{pro}
 \pf
 By the relation
 $
 \rho(x)(m)=\{\{x,-l\},m\}
 $
 and Lemma \ref{lem:relation}, (i) is obvious. (ii) follows from
\begin{eqnarray*}
 \rho(\xi)(m)&=&\{\{\xi,-l\},m\}=\{\{\xi,-l_1\},m\}=-\{\xi,\{l_1,m\}\}\\
 &=&\{\xi,\{m,l_1\}\}=\{\xi,l_1(m)\}=-\langle l_1^*(\xi),m\rangle.
\end{eqnarray*}
  (iii) follows from the relation
 $
 Dm=\{-l,m\}
 $
 and  Lemma \ref{lem:relation}.
(iv)-(vi) also follow from the relation
 $
\dorfman{a,b}=\{\{a,-l\},b\}
 $
 for all sections $a,b$ and Lemma \ref{lem:relation} . \qed\vspace{3mm}

 We have seen that we can obtain linear sections through the Dorfman bracket of constant sections. Thus, it is not enough if we  only consider constant sections.

 \begin{cor}
   Let $x\in\g_0$ be a constant section, and $\phi=\xi\otimes m,~\psi=\eta\otimes n:\g_{-1}^*\longrightarrow \g_0^*$ be linear sections, where $\xi,\eta\in\g_0^*$ and $m,n\in\g_{-1}$. Then we have
   \begin{itemize}
     \item[\rm(i)] $\dorfman{x,\xi\otimes m}=\xi\otimes l_2^1(x,m)+(l_2^0(x,\cdot)^* \xi)\otimes m;$
     \item[\rm(ii)] $\dorfman{\xi\otimes m,x}=-\xi\otimes l_2^1(x,m)-(l_2^0(x,\cdot)^* \xi)\otimes m+\xi(x)\big(l_1(m)-l_2^1(m,\cdot)\big);$
     \item[\rm(iii)] $\dorfman{\xi\otimes m,\eta\otimes n}=\langle l_1^*\eta,m\rangle\xi\otimes n-\langle l^*_1\xi,n\rangle\eta\otimes m$, equivalently, we have \begin{equation}\label{eq:special}\dorfman{\phi,\psi}=\phi\circ l_1^*\circ \psi-\psi\circ l_1^*\circ \phi.\end{equation}
   \end{itemize}
 \end{cor}
\pf Using the property of the Dorfman bracket, it is straightforward to obtain (i) and (ii). We only give the proof of (iii). We have
\begin{eqnarray*}
  \dorfman{\xi\otimes m,\eta\otimes n}&=&\eta\otimes \rho(\xi\otimes m)(n)+\dorfman{\xi\otimes m, \eta}\otimes n\\
  &=&\rho(\xi)(n)\eta\otimes m-\dorfman{\eta,\xi\otimes m}\otimes n\\
  &=&-l^*_1\xi(n)\eta\otimes m+l_1^*\eta(m)\xi\otimes n.
\end{eqnarray*}
Let it act on an arbitrary $\alpha\in \g_{-1}^*$. We obtain
$$
 \dorfman{\xi\otimes m,\eta\otimes n}(\alpha)=\alpha(n)l_1^*\eta(m)\xi-\alpha(m)l^*_1\xi(n)\eta=(\phi \circ l_1^* \circ \psi-\psi\circ l_1^*\circ \phi)(\alpha).
$$
This finishes the proof. \qed

\begin{rmk}
  Note that the bracket \eqref{eq:special} is naturally skew-symmetric and satisfies the Jacobi identity. Thus, it is a Lie bracket. It has already appeared in other places, e.g. see \cite[Proposition 3.1]{shengzhu2} for more details.
\end{rmk}
%From a Lie 2-algebra $(g,l)$, we get a Courant algebroid $E\longrightarrow{M}$, where
%$E=g^*_{-1}\times(g^*_0\oplus{g_0}), M=g^*_{-1}$.
Given  a Courant algebroid $E\longrightarrow{M}$, using the skew-symmetric Courant bracket, we get a $2$-term $L_\infty$-algebra structure on $C^\infty(M)\oplus{\Gamma(E)}$. Now consider the Courant algebroid \eqref{eq:E} obtained from a $2$-term $L_\infty$-algebra. Since it is linear, we
pick  linear functions on $\g_{-1}^*$ as the degree $-1$ part and $\mathfrak{g}_0\oplus{(\mathfrak{g}_{-1}}\otimes{\mathfrak{g}^*_0})$ as the degree $0$ part. Before we give the precise structure of the corresponding  2-term $L_\infty$-algebra, we give some properties of the skew-symmetric Courant bracket $\Courant{\cdot,\cdot}$.

\begin{cor}\label{cor:skewsym bracket}
  With the above notations, for any $x,y,z\in\g_0$, and $\xi\otimes m,\eta\otimes n\in\g_0^*\otimes \g_{-1}$, we have
\begin{equation}\label{eq:l2}
\left\{\begin{array}{rcl} \Courant{x,y}&=&x\circ y=l_2^0(x,y)+l_3(x,y,\cdot);\\
\Courant{x,\xi\otimes m}&=&\xi\otimes l_2^1(x,m)+(l_2^0(x,\cdot)^* \xi)\otimes m-\half\xi(x)\big(l_1(m)-l_2^1(m,\cdot)\big);\\
 \Courant{\xi\otimes m,\eta\otimes n}&=&l_1^*\eta(m)\xi\otimes n-l^*_1\xi(n)\eta\otimes m.\end{array}\right.
  \end{equation}
\end{cor}

Now we can give the main theorem in this section.
\begin{thm}\label{thm:newLie2}
Given a $2$-term $L_\infty$-algebra $\mathfrak{g}=(\g_{-1}\stackrel{l_1}{\longrightarrow}\g_0,l_2=l^0_2+l^1_2,l_3)$, we can obtain a new $2$-term $L_\infty$-algebra $\tilde{\g}=(\mathfrak{g}_{-1}\stackrel{\tilde{l}_1}{\longrightarrow}{\mathfrak{g}_0\oplus{(\mathfrak{g}_{-1}\otimes{\mathfrak{g}^*_0})}},
\tilde{l_2}=\tilde{l}_2^0+\tilde{l}_2^1,\tilde{l}_3)$ from the corresponding Courant algebroid \eqref{eq:E}, in which $\tilde{l}_1=D$ (given by \eqref{eq:D}), $\tilde{l}^0_2$ is given by \eqref{eq:l2}, $\tilde{l}^1_2$ and $\tilde{l}_3$ are given by
\begin{eqnarray}
\label{eq:bracket10}\tilde{l}^1_2(x+\xi\otimes m,n)&=&\frac{1}{2}l_2^1(x,n)+\half \langle\xi,l_1(n)\rangle m,
\end{eqnarray}
and
\begin{eqnarray}
\nonumber&&\tilde{l}_3(x_1+\xi_1\otimes m_1,x_2+\xi_2\otimes m_2,x_3+\xi_3\otimes m_3)\\
\nonumber&=&-\half l_3(x_1,x_2,x_3)-\half\big(\langle l_2^0(x_1,x_2),\xi_3\rangle m_3+c.p.\big)\\
\nonumber&&-\frac{1}{4}\big(\langle\xi_1,x_2\rangle\langle\xi_3,l_1m_1\rangle m_3-\langle\xi_2,x_1\rangle\langle\xi_3,l_1m_2\rangle m_3+c.p.\big)\\
\label{eq:bracketl3}&&-\frac{1}{4}\big( \langle\xi_2,x_3\rangle l_2^1(x_1,m_2)-\langle\xi_3,x_2\rangle l_2^1(x_1,m_3)+c.p.\big).
\end{eqnarray}
\end{thm}
\pf  It follows from Theorem \ref{thm:CALie2}, Proposition \ref{pro:CA structure} and Corollary \ref{cor:skewsym bracket}. \qed\vspace{3mm}

Now we give some examples to illustrate that the construction given in the above theorem is full of interest.

\begin{ex}\label{ex:1}\rm{({\bf omni-Lie algebra}) Let $V$ be a vector space. Consider the abelian $2$-term $L_\infty$-algebra
 $(V\stackrel{\Id}{\longrightarrow} V,l_2=0,l_3=0)$. By Theorem \ref{thm:newLie2}, we get a new $2$-term $L_\infty$-algebra $(V\stackrel{\id}\longrightarrow V\oplus \gl(V),\tilde{l}_2,\tilde{l}_3)$, in which $\id$ is the natural inclusion, and
\begin{equation*}
  \left\{\begin{array}{rcll}
 \tilde{l}_2^1(u+A,m)&=&\frac{1}{2}Am,\\
 \tilde{l}_2^0(u+A,v+B)&=&\frac{1}{2}(Av-Bu)+[A,B],\\
  \tilde{l}_3(u+A,v+B,w+C)&=&-\frac{1}{4}([A,B]w+[B,C]u+[C,A]v),
   \end{array}\right.
\end{equation*}
for all $u,v,w\in V_0=V, m\in V_{-1}=V$ and $A,B,C\in \gl(V)$.
This $2$-term $L_\infty$-algebra is the one associated to the omni-Lie algebra $V\oplus \gl(V)$. See \cite{shengzhu1,Alan} for more details.
}
\end{ex}

\begin{ex}{\rm({\bf 2-term $L_\infty$-algebra of string type})
Let  $(\mathfrak{k},[\cdot,\cdot]_\frkk)$ be a Lie algebra. Consider the 2-term $L_\infty$-algebra $(\mathbb{R}\stackrel{0}\longrightarrow \mathfrak{k},l_2=[\cdot,\cdot]_\frkk,l_3=0)$. By Theorem \ref{thm:newLie2}, we get a new $2$-term $L_\infty$-algebra $(\mathbb{R}\stackrel{0}\longrightarrow \mathfrak{k}\oplus\mathfrak{k}^*,\tilde{l}_2,\tilde{l}_3)$, where $\tilde{l}_2$ and $\tilde{l}_3$ are given by
\begin{equation*}
  \left\{\begin{array}{rcll}
 \tilde{l}_2^1(u+\xi,r)&=&0,\\
 \tilde{l}_2^0(u+\xi,v+\eta)&=&[u,v]_\frkk+ \ad_u^*\eta-\ad_v^*\xi,\\
  \tilde{l}_3(u+\xi,v+\eta,w+\zeta)&=&-\frac{1}{2}(\langle[u,v]_\frkk,\zeta\rangle
  +\langle[v,w]_\frkk,\xi\rangle+\langle[w,u]_\frkk,\eta\rangle),
   \end{array}\right.
\end{equation*}
for all $u,v,w\in\frkk$, $\xi,\eta,\zeta\in\frkk^*$ and $r\in\mathbb R$. This is exactly the 2-term $L_\infty$-algebra of string type studied in
 \cite{shengzhu1}.}
\end{ex}

\begin{ex}{\rm
Let $(\mathfrak{k},[\cdot,\cdot]_\frkk,K)$ be a quadratic Lie algebra, and $(\mathbb{R}\stackrel{0}\longrightarrow \mathfrak{k},l_2,l_3)$ the corresponding string Lie 2-algebra. More precisely, $$l_2^0(u,v)=[u,v]_\frkk,\  l_2^1(u,r)=0,\
l_3(u,v,w)=K([u,v]_\frkk,w), \ \forall u,v,w\in \mathfrak{k},\ r\in \mathbb{R}.$$ Then we get a new $2$-term $L_\infty$-algebra $(\mathbb{R}\stackrel{0}\longrightarrow \mathfrak{k}\oplus\mathfrak{k}^*,\tilde{l}_2,\tilde{l}_3)$, in which
\begin{equation*}
  \left\{\begin{array}{rcll}
 \tilde{l}_2^1(u+\xi,r)&=&0,\\
 \tilde{l}_2^0(u+\xi,v+\eta)&=&[u,v]_\frkk+ K([u,v]_\frkk,\cdot)+\ad_u^*\eta-\ad_v^*\xi,\\
  \tilde{l}_3(u+\xi,v+\eta,w+\zeta)&=&-\frac{1}{2}K([u,v]_\frkk,w)-\frac{1}{2}(\langle[u,v]_\frkk,\zeta\rangle
  +\langle[v,w]_\frkk,\xi\rangle+\langle[w,u]_\frkk,\eta\rangle),
   \end{array}\right.
\end{equation*}
for all $u,v,w\in \mathfrak{k}, r\in \mathbb{R}$ and $\xi,\eta,\zeta\in \mathfrak{k}^*$.
}\end{ex}
\begin{ex}\label{ex:3}\rm{
Consider  the $2$-term $L_\infty$-algebra $(\mathfrak{k}\stackrel{\Id}\longrightarrow \mathfrak{k},l_2=[\cdot,\cdot]_\frkk,l_3=0)$, where $(\mathfrak{k},[\cdot,\cdot]_\frkk)$ is a Lie algebra. By Theorem \ref{thm:newLie2}, we obtain a new $2$-term $L_\infty$-algebra $(\mathfrak{k}\stackrel{\tilde{l}_1}\longrightarrow \mathfrak{k}\oplus \gl(\mathfrak{k}),\tilde{l}_2,\tilde{l}_3)$, where
\begin{equation*}
  \left\{\begin{array}{rcll}
\tilde{l}_1(m)&=&m-\ad_m,\\
 \tilde{l}_2^1(u+A,m)&=&\frac{1}{2}[u,m]_\frkk+\frac{1}{2}Am,\\
 \tilde{l}_2^0(u+A,v+B)&=&[u,v]_\frkk+\frac{1}{2}(Av-Bu)+[\ad_u, B]+[ A,\ad_v]+\frac{1}{2}(\ad_{Bu}-\ad_{Av})+[A,B],\\
  \tilde{l}_3(u+A,v+B,w+C)&=&-\frac{1}{2}C[u,v]_\frkk-\frac{1}{4}[A,B]w-\frac{1}{4}([u,Bw]_\frkk+[Bu,w]_\frkk)+c.p.,
   \end{array}\right.
\end{equation*}
for all $u,v,w\in \frkk, m\in \frkk$ and $A,B,C\in \gl(\frkk)$.
This 2-term $L_\infty$-algebra can be viewed as a deformation of the one associated to the omni-Lie algebra. It can also be viewed as a linearization of the Courant algebroid $TM\oplus T^*_\pi M$, where $T^*_\pi M$ denotes the Lie algebroid associated to the Poisson manifold $(M,\pi)$.  We will study its full properties in a separate paper.}
\end{ex}

Now let us consider the relation between the $2$-term $L_\infty$-algebra given in the above theorem and the original $2$-term $L_\infty$-algebra.
\begin{thm}\label{thm:morphism}
With the above notations, there is an $L_\infty$-algebra homomorphism $F=(F_0=\pr_1,F_1={\Id},F_2)$ from $\tilde{\g}$ to $\g$, in which $\pr_1$ is the projection to the first component, and
$F_2:\wedge^2(\g_0\oplus(\g_0^*\otimes \g_{-1}))\longrightarrow\g_{-1}$ is given by
 \begin{equation}
F_2(x+\xi\otimes m,y+\eta\otimes n)=\half\big(\langle{\xi,y}\rangle m-\langle{x,\eta}\rangle n\big)=\half\pairm{x+\xi\otimes m,y+\eta\otimes n},
 \end{equation}
 for all $ x,y\in{\mathfrak{g}_0},~\xi\otimes m,\eta\otimes n\in\g_0^*\otimes \g_{-1}.$
\end{thm}
\pf We verify the conditions of an $L_\infty$-algebra homomorphism step by step. It is obvious that
\begin{equation}\label{eq:mor0}
F_{0}\circ{D}=l_{1}\circ{F_1}.
\end{equation}
Then,
by Corollary \ref{cor:skewsym bracket}, we obtain
\begin{eqnarray*}
F_0\tilde{l}_2^0(x+\xi\otimes m,y+\eta\otimes n)&=&l_2^0(x,y)+\frac{1}{2}\big(\langle{\xi,y}\rangle l_{1}m-\langle{x,\eta}\rangle l_{1}n\big).
\end{eqnarray*}
On the other hand, it is obvious that
$$
l_2^0(F_0( x+\xi\otimes m),F_0(y+\eta\otimes n))=l_2^0(x,y).
$$
Therefore, we have
\begin{eqnarray}
 \nonumber &&F_0\tilde{l}_2^0(x+\xi\otimes m,y+\eta\otimes n)-l_2^0(F_0( x+\xi\otimes m),F_0(y+\eta\otimes n))\\\nonumber &=&\frac{1}{2}\big(\langle{\xi,y}\rangle l_{1}m-\langle{x,\eta}\rangle l_{1}n\big)\\
 \label{eq:mor1}&=&l_1 F_2(x+\xi\otimes m,y+\eta\otimes n).
\end{eqnarray}

By \eqref{eq:bracket10}, we obtain
$$
F_1\tilde{l}^1_2(x+\xi\otimes m,n)=\frac{1}{2}l_2^1(x,n)+\half\langle\xi,l_1(n)\rangle m,
$$
which implies that
$$
 F_1\tilde{l}_2^1(x+\xi\otimes m, n)-l_2^1(F_0( x+\xi\otimes m),F_1(n))=-\frac{1}{2}l_2^1(x,n)+\half \langle\xi,l_1(n)\rangle m.
$$
On the other hand, we have
$$
F_2(x+\xi\otimes m, Dn)=F_2(x+\xi\otimes m, l_1(n)-l_2^1(n,\cdot))=-\frac{1}{2}l_2^1(x,n)+\half \langle\xi,l_1(n)\rangle m.
$$
Thus, we have
\begin{equation}\label{eq:mor2}
 F_1\tilde{l}_2^1(x+\xi\otimes m, n)-l_2^1(F_0( x+\xi\otimes m),F_1(n))=F_2(x+\xi\otimes m, Dn).
\end{equation}

It is left to prove the compatibility condition, that is, for all $a,b,c\in\g_0\oplus\g_0^*\otimes \g_{-1}$,
\begin{equation}\label{eq:mor3}
l_2^1(F_{0}a,F_2(b,c))+c.p.+l_3(F_{0}a,F_{0}b,F_{0}c)=F_2(\tilde{l}_2^0(a,b),c)+c.p.+F_1{\tilde{l}_3}(a,b,c).
\end{equation}
First, for all $a,b,c\in\mathfrak{g}_{-1}\otimes{\mathfrak{g}^*_{0}}$, it is obvious that both sides of the equation are $0$. While for all $a,b,c\in{\mathfrak{g}_0}$ (we use the notation $x,y,z$ instead), it is obvious that the left hand side is equal to $l_3(x,y,z)$. It is not hard to see that
\begin{eqnarray*}
F_2(\tilde{l}_2^0(x,y),z)+c.p.&=&F_2(l_3(x,y,\cdot),z)+c.p.=\frac{1}{2}l_3(x,y,z)+c.p.=\frac{3}{2}l_3(x,y,z),\\
F_1{\tilde{l}_3}(x,y,z)&=&=-\frac{1}{2}l_3(x,y,z),\quad \mbox{by ~ Theorem ~\ref{thm:newLie2}}
\end{eqnarray*}
which implies that \eqref{eq:mor3} holds for all $x,y,z\in{\mathfrak{g}_0}$.
For $x,y$ in $\mathfrak{g}_0$ and  $\xi\otimes m$ in $\g_0^*\otimes \mathfrak{g}_{-1}$, we have
\begin{eqnarray*}
l_2^1(F_{0}x,F_2(y,\xi\otimes m))+c.p.+l_3(F_{0}x,F_{0}y,F_{0}(\xi\otimes m))&=&l_2^1(x,-\frac{1}{2}\langle{\xi,y}\rangle m)+l_2^1(y,\frac{1}{2}\langle{\xi,x}\rangle m)\\&=&
-\frac{1}{2}\langle{y,\xi}\rangle{l_2^1}(x,m)+\frac{1}{2}\langle{x,\xi}\rangle{l_2^1}(y,m),
\end{eqnarray*}
and
\begin{eqnarray*}
&&F_2(\tilde{l}_2^0(x,y),\xi\otimes m)+c.p.+F_1{\tilde{l}_3}(x,y,\xi\otimes m)\\&=&\big(-\frac{1}{2}-\frac{1}{6}\big)\langle{\Courant{x,y},\xi}\rangle m+\big(\frac{1}{2}-\frac{1}{6}\big)\langle{\Courant{y,\xi\otimes m},x}\rangle+
\big(\frac{1}{2}-\frac{1}{6}\big)\langle{\Courant{\xi\otimes m,x},y}\rangle\\&=&-\frac{2}{3}\langle{l_2^0(x,y),\xi}\rangle m+\frac{1}{3}\big\{
\langle{\xi,x}\rangle l_2^1(y,m)-\langle{\xi,l_2^0(y,x)}\rangle m+\frac{1}{2}\langle{\xi,y}\rangle l_2^1(m,x)\big\}\\&&+
\frac{1}{3}\big\{
-\langle{\xi,y}\rangle l_2^1(x,m)+\langle{\xi,l_2^0(x,y)}\rangle m-\frac{1}{2}\langle{\xi,x}\rangle l_2^1(m,y)\big\}\\&=&
-\frac{1}{2}\langle{y,\xi}\rangle{l_2^1}(x,m)+\frac{1}{2}\langle{x,\xi}\rangle{l_2^1}(y,m),
\end{eqnarray*}
which implies that \eqref{eq:mor3} holds for two elements in $\g_0$ and one element in $\g_0^*\otimes \mathfrak{g}_{-1}$. At last, for  $x\in\mathfrak{g}_0$  and $\xi\otimes m,\eta\otimes n\in\g_0^*\otimes \mathfrak{g}_{-1}$, it is obvious that the left hand side of \eqref{eq:mor3} is equal to $0$. Furthermore, we have
\begin{eqnarray*}
&&F_2(\tilde{l}_2^0(x,\xi\otimes m),\eta\otimes n)+c.p.+F_1{\tilde{l}_3}(x,\xi\otimes m,\eta\otimes n)\\&=&\big(-\frac{1}{2}-\frac{1}{6}\big)\langle{\Courant{x,\xi\otimes m},\eta}\rangle n+\big(\frac{1}{2}
-\frac{1}{6}\big)\langle{\Courant{\xi\otimes m,\eta\otimes n},x}\rangle+\big(-\frac{1}{2}
-\frac{1}{6}\big)\langle{\Courant{\eta\otimes n,x},\xi}\rangle m\\&=&-\frac{2}{3}\times \big(-\frac{1}{2}\big)\langle{\xi,x}\rangle \langle{l_1 m,\eta}\rangle n+\frac{1}{3}\big\{-\langle{l_1 m,\eta}\rangle \langle{\xi ,x}\rangle n +\langle{l_1 n,\xi}\rangle \langle{\eta,x}\rangle m\big\}\\ &&
-\frac{2}{3}\times \frac{1}{2}\langle{l_1 n,\xi}\rangle \langle{\eta,x}\rangle m\\&=&0.
\end{eqnarray*}
Thus, \eqref{eq:mor3} holds for all $a,b,c\in\g_0\oplus\g_0^*\otimes \g_{-1}$. By \eqref{eq:mor0}-\eqref{eq:mor3}, we deduce that $F$ is an $L_\infty$-algebra homomorphism. \qed

\section{Lie $2$-algebras and quasi-Poisson groupoids}

In this section, we show that the Courant algebroid associated to a 2-term $L_\infty$-algebra is the double of a Lie quasi-bialgebroid. Then, by integration, we obtain a quasi-Poisson groupoid. Since a 2-term $L_\infty$-algebra is equivalent to a Lie 2-algebra, this quasi-Poisson groupoid can serve as the geometric structure on the dual of a Lie 2-algebra.

 Let $\mathfrak{g}=(\g_{-1}\stackrel{l_1}{\longrightarrow}\g_0,l_2=l^0_2+l^1_2,l_3)$ be a $2$-term $L_\infty$-algebra. The vector bundle $E\longrightarrow{M}$ given by \eqref{eq:E} can be decomposed as $A\oplus{A^*}$, where $A={\mathfrak{g}^*_{-1}}\times\mathfrak{g}^*_{0}\longrightarrow{\mathfrak{g}^*_{-1}}$ is a trivial vector bundle. By Lemma \ref{lem:delta-big-bracket} ($k=-1$), the $2$-term $L_\infty$-algebra structure $l=l_1+(l^1_2+l^0_2)+l_3$, where  $l_1\in{\mathfrak{g}^*_{-1}\otimes{\mathfrak{g}_0}},~ l^0_2
\in{\wedge^{2}\mathfrak{g}^*_{0}\otimes{\mathfrak{g}_0}},~ l^1_2\in{\mathfrak{g}^*_{0}\wedge{\mathfrak{g}^*_{-1}}\otimes{\mathfrak{g}_{-1}}},~ l_3
\in{\wedge^{3}\mathfrak{g}^*_{0}\otimes{\mathfrak{g}_{-1}}},$ satisfies $\{l,l\}=0$. By the degree reason, $\{l,l\}=0$ implies that
\begin{eqnarray}
\label{Eqn:p L d mu k}
\{l_1,l_1\}=0\quad
\{l_1,l_2\}=0,\quad
\half\{l_2,l_2\}+\{l_1,l_3\}=0,\quad
\{l_2,l_3\}=0.
\end{eqnarray}
Comparing this with \eqref{eq:defiquasi}, we can see that Eq. \eqref{Eqn:p L d mu k} implies Eq. \eqref{eq:defiquasi} for $\mu=-l_1,~\gamma=-l_2=-l_2^0-l_2^1,~\phi=-l_3$. Therefore, we have

\begin{thm}\label{thm:Lie quasi}
  Let $\mathfrak{g}=(\g_{-1}\stackrel{l_1}{\longrightarrow}\g_0,l_2=l^0_2+l^1_2,l_3)$ be a $2$-term $L_\infty$-algebra. Then $(A,-l)$ is a Lie quasi-bialgebroid, where $A={\mathfrak{g}^*_{-1}}\times\mathfrak{g}^*_{0}\longrightarrow{\mathfrak{g}^*_{-1}}$ and $l=l_1+(l^1_2+l^0_2)+l_3$, and the Courant algebroid $E$ is the double of the Lie quasi-bialgebroid $(A,-l)$.
\end{thm}

Now we give the precise structure on $A$ and $A^*$ using the general method of differential geometry.
By Theorem \ref{thm:Lie quasi} and Proposition \ref{pro:CA structure}, we have

\begin{cor}\label{cor:Lie algebroid}
  The Lie algebroid structure on $A$, determined by $-l_1$, is given by
  \begin{itemize}
  \item[\rm(i)] for any constant section $\xi\in\g_0^*,$  the anchor $\rho_A$ is given by $\rho_A(\xi)=-l_1^*(\xi)$;
    \item[\rm(ii)] for any constant sections $\xi,\eta\in\g_0^*,$ we have
    $[\xi,\eta]_A=0;$
    \item[\rm(iii)] for any constant section $\xi\in\g_0^*$ and linear section $\eta\otimes n\in\g_0^*\otimes\g_{-1}$, we have
    $$[\xi,\eta\otimes n]_A=\langle\xi,l_1(n)\rangle\eta;$$
    \item[\rm(iv)] for any linear sections $\xi\otimes m, \eta\otimes n\in\g_0^*\otimes\g_{-1}$, we have
    $$[\xi\otimes m,\eta\otimes n]_A=\langle l_1^*\eta,m\rangle\xi\otimes n-\langle l^*_1\xi,n\rangle\eta\otimes m.$$
  \end{itemize}
  Thus, $A$ is an action Lie algebroid of the abelian Lie algebra $\g_0^*$ acting on $\g_{-1}^*$ via $-l^*_1$, which sends an element $\xi\in\g_0^*$ to a constant vector field $-l^*_1(\xi)\in\g_{-1}^*.$ The corresponding Chevalley-Eilenberg coboundary operator $d_A:\Gamma(\wedge^\bullet A^*)\longrightarrow\Gamma(\wedge^{\bullet+1}A^*)$ is determined by
  $$
  d_Am=\{-l_1,m\}=l_1(m).
  $$
\end{cor}
\begin{cor}\label{cor:Lie quasialgebroid}
  For all constant section $x\in\g_0$ of $A^*$, $-l_2^1$ gives rise to the anchor map $\rho_{A^*}$ of $A^*$ via
  $$\rho_{A^*}(x)=l_2^1(x,\cdot),$$
  which is a linear vector field. For all constant sections $x,y\in\g_0$, $-l_2^0$ gives rise to the bracket operation on $A^*$:
  $$
  [x,y]_{A^*}=l_2^0(x,y).
  $$
 The Jacobi identity of $[\cdot,\cdot]_{A^*}$ is controlled by $\phi=-l_3\in\wedge^3\g_0^*\otimes \g_{-1}\subset\Gamma(\wedge^3A)$. More precisely, we have
  $$
  [[x,y]_{A^*},z]_{A^*}+c.p.=d_A\phi(x,y,z)+\phi(d_A x,y,z)-\phi(x,d_A y,z)+\phi(x,y,d_A z).
  $$
\end{cor}
Now we can restate Theorem \ref{thm:Lie quasi} using usual language of differential geometry. This is necessary since we will consider the integration of Lie quasi-bialgebroids at the end of this section.

\begin{cor}
Given a $2$-term $L_\infty$-algebra $\mathfrak{g}=(\g_{-1}\stackrel{l_1}{\longrightarrow}\g_0,l_2=l^0_2+l^1_2,l_3)$, we obtain a Lie quasi-bialgebroid $(A,\delta,\phi)$, where the Lie algebroid $A=\mathfrak{g}^*_{-1}\times{\mathfrak{g}^*_{0}}\longrightarrow{\mathfrak{g}^*_{-1}}$ is given by Corollary \ref{cor:Lie algebroid}, $\delta:\Gamma(\wedge^kA)\longrightarrow \Gamma(\wedge^{k+1}A)$ is the generalized Chevalley-Eilenberg operator determined by the anchor $\rho_{A^*}$ and the bracket $[\cdot,\cdot]_{A^*}$ given in Corollary \ref{cor:Lie quasialgebroid}, and $\phi=-l_3$.
\end{cor}

In the last section, we have seen that, given a $2$-term $L_\infty$-algebra, there is a new $2$-term $L_\infty$-algebra associated to the Courant algebroid $E$ (Theorem \ref{thm:newLie2}), and there is an $L_\infty$-algebra homomorphism from the new $2$-term $L_\infty$-algebra to the original one (Theorem \ref{thm:morphism}). In this section, we see that the Courant algebroid is exactly the double of the Lie quasi-bialgebroid $(A,\delta,\phi)$. In fact, we can obtain a more general result than Theorem \ref{thm:morphism}. Namely, given any Lie quasi-bialgebroid $(A,\delta,\phi)$, we can obtain a $2$-term $L_\infty$-algebra $\CWM\oplus \Ker(d_A|_{\Gamma(A^*)})$, where $d_A$ is the coboundary operator associated to Lie algebroid $A$ and $d_A|_{\Gamma(A^*)}:\Gamma(A^*)\longrightarrow\Gamma(\wedge^2A^*)$ is the restriction of $d_A$ on $\Gamma(A^*)$. Furthermore, we construct an $L_\infty$-algebra homomorphism from the $2$-term $L_\infty$-algebra associated the Courant algebroid $A\oplus A^*$ to $\CWM\oplus \Ker(d_A|_{\Gamma(A^*)})$.

\begin{thm}\label{thm:dual L_inf}
Let $(A,\delta,\phi)$ be a Lie quasi-bialgebroid. Then there is a $2$-term $L_\infty$-algebra structure on  $\CWM\oplus \Ker(d_A|_{\Gamma(A^*)})$, where $\CWM$ is of degree $-1$, $\Ker(d_A|_{\Gamma(A^*)})$ is of degree $0$,
  and  $\{ l_i\}$ are given by
  \begin{equation}
  \left\{\begin{array}{rcll}
  l_1(f)&=&d_Af,&~\forall ~f\in\CWM,\\
   l_2(e_1,e_2)&=& [e_1,e_2]_{A^*} & ~\forall~e_1,e_2\in \Ker(d_A|_{\Gamma(A^*)}),\\
  l_2(e_1,f)&=& \rho_{A^*}(e_1)(f) &
  ~\forall~e_1\in \Ker(d_A|_{\Gamma(A^*)}),f\in\CWM,\\
  l_3(e_1,e_2,e_3)&=&-\phi(e_1,e_2,e_3) &~\forall~e_1,e_2,e_3\in\Ker(d_A|_{\Gamma(A^*)}).
   \end{array}\right.
\end{equation}

Furthermore, there is an $L_\infty$-algebra homomorphism $F=(F_0=\pr_2,F_1={\Id},F_2)$ from the $2$-term $L_\infty$-algebra associated to the Courant algebroid $A\oplus A^*$, $\CWM\oplus \Big(\Gamma(A)\oplus\Ker(d_A|_{\Gamma(A^*)})\Big)$, to $\CWM\oplus \Ker(d_A|_{\Gamma(A^*)})$, where $\pr_2$ is the projection to the second component, and $F_2$ is given by
\begin{equation}
  F_2(e_1,e_2)=\half\pairm{e_1,e_2}, \quad\forall e_1,e_2\in \Gamma(A)\oplus\Ker(d_A|_{\Gamma(A^*)}).
\end{equation}
\end{thm}
\pf Since $A$ is a Lie algebroid, we have $d_A^2=0,$ i.e. $\Img(d_A)\subset \ker(d_A)$. Thus, $l_1$ is well-defined. By the fact that $(A,\delta,\phi)$ is a Lie quasi-bialgebroid, it is straightforward to see that $\{l_i\}$ is a $2$-term $L_\infty$-algebra structure. Furthermore, it is not hard to see that the Courant bracket on $\Gamma(A)\oplus\Ker(d_A|_{\Gamma(A^*)})$ is closed, which implies that the $2$-term $L_\infty$-algebra $\CWM\oplus\Big(\Gamma(A)\oplus\Gamma(A^*)\Big)$ associated to the Courant algebroid $A\oplus A^*$ can be reduced to $\CWM\oplus \Big(\Gamma(A)\oplus\Ker(d_A|_{\Gamma(A^*)})\Big)$. We can also prove that $F$ is indeed a homomorphism similar as the proof of Theorem \ref{thm:morphism}. \qed\vspace{3mm}

Next, we return to the Lie quasi-bialgebroid $(A,\delta,\phi)$ associated to a $2$-term $L_\infty$-algebra. Clearly, the Lie algebroid $A$ can be integrated to an action groupoid $\Gamma:\mathfrak{g}^*_{-1}\times{\mathfrak{g}^*_0}\rightrightarrows{\mathfrak{g}^*_{-1}}$ with the abelian group structure on $\mathfrak{g}^*_0$, where the source, target and inclusion maps are given by
$$s(\alpha,\xi)=\alpha, ~t(\alpha,\xi)=\alpha+l^*_{1}\xi, ~ i(\alpha)=(\alpha,0),$$ for all $(\alpha,\xi)\in \mathfrak{g}^*_{-1}\times{\mathfrak{g}^*_0}$. By Theorem $4.9$ in \cite{Xu},
we know that $\Gamma$ has a quasi-Poisson structure, such that its corresponding Lie quasi-bialgebroid is exactly $(A,\delta,\phi)$. Recall that a quasi-Poisson groupoid is a triple $(\Gamma,\Pi,\phi)$, where $\Gamma$ is a Lie groupoid whose Lie algebroid is $A$, $\Pi\in \wedge^2\frkX(\Gamma)$, and $\phi\in\Gamma(\wedge^3A)$, satisfying   $\frac{1}{2}[\Pi,\Pi]_S=\overleftarrow{\phi}-\overrightarrow{\phi}, [\Pi,\overleftarrow{\phi}]=0$.
Generally, it is a rather difficult work to get the bivector field from the data of an ordinary Lie quasi-bialgebroid.
Nevertheless, we can elaborate the quasi-Poisson structure in our specific case, since all the
structures are determined by the information of the $2$-term $L_\infty$-algebra.

\begin{thm}\label{thm:int}
The quasi-Poisson groupoid corresponding to $(A,\delta,\phi)$ is $(\Gamma,\Pi, \phi=-l_3)$, where $\Pi$ is characterized by
\[\Pi(dx,dy)=-l_2^0(x,y), \ \ \Pi(dx,dm)=-l_2^1(x,m), \ \ \Pi(dm,dn)=-l_2^1(l_1(m),n),\]
%\comment{I drop the minus sigh in the first equation, or $\delta\xi$ differs a minus sigh from $[\Pi,\xi]_S$ in the last part of the proof. If you think this is not natural, contact me.}
where $d$ is the usual de Rham differential, and $x,y\in \g_0, m,n\in \g_{-1}$ are  linear functions on
$\mathfrak{g}^*_{-1}\times{\mathfrak{g}^*_0}$.
\end{thm}
We need to make some preparations to prove the theorem. Let us take a general action groupoid  into account.

Let $M{\lhd}G$ be an action groupoid and $M{\lhd}\g$ the corresponding Lie algebroid, where
$\g$ is the Lie algebra of Lie group $G$, and the anchor is given by the infinitesimal action of $\g$ on  $M$, i.e. a Lie algebra homomorphism $\widehat{\cdot}:{\g}\longrightarrow{\frkX(M)}$,
\begin{equation*}
\xymatrix@C=0.5cm{  \widehat{X}(x)=\frac{d}{dt}x\cdot{\exp(-tX)}, \quad \quad  \forall X\in{\g},x\in{M}
.                }
\end{equation*}
We use $l_a(resp. ~r_a)$ to denote the left(right)-translation on the Lie group $G$ and $l_{(x,a)}(resp. ~r_{(x,a)})$ the left(right)-translation on the groupoid $M\triangleleft G$ respectively. Let $(l_a)_*,(r_a)_*$ and  $(l_{(x,a)})_*,(r_{(x,a)})_*$ be the corresponding tangent maps.

For the Lie algebroid $A$ of a Lie groupoid $\Gamma$, when we say the left-translation $\overleftarrow{\Lambda}$ (right-translation $\overrightarrow{\Lambda}$) of ${\Lambda}\in \Gamma(A)$, we are identifying $\Gamma(A)$ with left-invariant (right-invariant) vector fields on $\Gamma$. We assume that the derivatives are all taken at $t=0$.

\begin{lem}\label{lem:ltrt}
Let $M{\lhd}G$ be an action Lie groupoid and $M{\lhd}\g$ the corresponding Lie algebroid. For all ${\phi}\in{\Gamma(M{\lhd}\g)},(x,a)\in{M\times{G}}$, we have

$$\overleftarrow{\phi}_{(x,a)}=(0_{x},(l_a)_{*e}\phi_{x\cdot{a}}), \quad \quad
\overrightarrow{\phi}_{(x,a)}=((\widehat{\phi_x})_x,(r_a)_{*e}\phi_x).
$$
\end{lem}
\pf For this action groupoid, the $s$-fiber over $x$ is $s^{-1}(x)=\{(x,a);\forall{a}\in{G}\}=\{x\}\times{G}$, and the $t$-fiber over $x$ is
$t^{-1}(x)=\{(x\cdot{a^{-1}},a);\forall{a}\in{G}\}$.
 The tangent map of $l_{(x,a)}:s^{-1}(x\cdot{a})\longrightarrow{s^{-1}(x)}$ is  $(l_{(x,a)})_{*(x\cdot a,e)}:T_{(x\cdot{a},e)}s^{-1}(x\cdot{a})
\longrightarrow{T_{(x,a)}s^{-1}(x)}$.
Let $\phi:M\longrightarrow \g$ be a section of $M\lhd\g$ and $\gamma(t)$ a smooth curve on $G$ satisfying  $\gamma(0)=e,\gamma'(0)=\phi(x\cdot{a})$. Then we have
\begin{eqnarray*}
\overleftarrow{\phi}_{(x,a)}=(l_{(x,a)})_* \phi_{x\cdot{a}}=\frac{d}{dt}(x,a)(x\cdot{a},\gamma(t))=\frac{d}{dt}(x,a\gamma(t))
=(0_{x},(l_a)_{*e}\phi_{x\cdot{a}}).
\end{eqnarray*}

Similarly, the tangent map of $r_{(x,a)}:t^{-1}(x)\longrightarrow{t^{-1}(x\cdot{a})}$ is
$(r_{(x,a)})_{*(x,e)}:T_{(x,e)}t^{-1}(x)\longrightarrow{T_{(x,a)}t^{-1}(x\cdot{a})}$. Let $\gamma(t)$ be a smooth curve on $G$ satisfying $\gamma(0)=e,\gamma'(0)=\phi(x)$. When $\phi$ is considered as a right-invariant vector field, we have
\[\phi_{x}=\frac{d}{dt}(x,\gamma(-t))^{-1}=\frac{d}{dt}(x\cdot{\gamma(-t)},\gamma(-t)^{-1}).\]
Thus, we get
\begin{eqnarray*}
\overrightarrow{\phi}_{(x,a)}&=&(r_{(x,a)})_* \phi_{x}=\frac{d}{dt}(x\cdot{\gamma(-t)},\gamma(-t)^{-1})(x,a)\\&=&\frac{d}{dt}(x\cdot{\gamma(-t)},
\gamma(t)a)=((\widehat{\phi_x})_x,(r_a)_{*e}\phi_x).\qed\vspace{3mm}
\end{eqnarray*}

$\mathbf{Proof \ of \ Theorem\ \ref{thm:int}:}$ According to Theorem 2.34 and Theorem 4.9 in \cite{Xu}, we only need to prove $\delta_\Pi=\delta,$ that is,
for any function $f\in C^\infty(\mathfrak{g}_{-1}^*)$  and  any section $e\in \Gamma(A) $, we have
 $\overleftarrow{\delta f}=[t^* f,\Pi]_S,~ \overleftarrow{\delta e}=[\overleftarrow{e},\Pi]_S$. We only give the proof for
a linear function $m\in \mathfrak{g}_{-1}$  and a constant section $\xi\in \mathfrak{g}_0^* $.

 By calculation, we have $t^* (m)=m-l_1(m).$ Considering the linear action of the abelian group $\g_0^*$ on $\g_ {-1}^*$ and by Lemma \ref{lem:ltrt}, we have,
\begin{eqnarray*}
  \overleftarrow{\delta m}&=&\delta m -l_1 (\delta m)=-l_2^1(m,\cdot)+l_1 (l_2^1(m,\cdot)),\\
  \overleftarrow{\delta \xi}&=&\delta \xi -l_1 (\delta \xi)=\delta \xi,\\
  \overleftarrow{\xi}&=&\xi.
\end{eqnarray*}
 Thus, for a linear function $n+x$ on $\g_{-1}^* \times \g_0^*$, we have
 $$
\overleftarrow{\delta m}(n+x)=-l_2^1(m,x)+l_1 (l_2^1(m,x)),
$$
while
\begin{eqnarray*}
[t^* m,\Pi]_S(n+x)&=&\Pi(dm,dn)+\Pi(dm,dx)-\Pi(d l_1(m),dn)-\Pi(d l_1 (m),dx)\\&=&-l_2^1(m,x)+l_2^0(l_1(m),x).
\end{eqnarray*}
Then, we have $\overleftarrow{\delta m}=[t^* m,\Pi]_S$.
To prove $\delta \xi=[\xi,\Pi]_S$ as bivector fields, we just need to verify that it holds on two 1-forms on $\g_0^*$, since it vanishes for all the other cases. By straightforward calculation, we have
$
\delta \xi(dx,dy)=-\langle\xi,[x,y]_{A^*}\rangle=-\langle\xi,l_2^0(x,y)\rangle,$ and
\begin{eqnarray*}
[\xi,\Pi]_S(dx,dy)&=&L_\xi \Pi(dx,dy)\\&=&-\langle\xi,l_2^0(x,y)\rangle-\Pi(L_\xi dx,dy)-\Pi(dx,L_\xi dy)
%\\&=&\langle\xi,l_2^0(x,y)\rangle-2\langle\xi,l_2^0(x,y)\rangle+0
\\&=&-\langle\xi,l_2^0(x,y)\rangle.
\end{eqnarray*}
%\comment{Here, $[\Pi,\xi]=-[\xi,\Pi]$, $\Pi(L_\xi dx,dy)+\Pi(dx,L_\xi dy)=0$. Thus the two results differ a minus sigh. One way is to }
Therefore, we complete the proof.\qed\vspace{3mm}

There is a one-to-one correspondence between  2-term chain complex ${V}_{1}\stackrel{\dM}{\longrightarrow}{V}_{0}$ and 2-vector space $\mathbb{V}:{V}_{1}\oplus{V_0}\rightrightarrows{V_0}$. Actually, the 2-vector space $\mathbb{V}$ is an action Lie groupoid $V_0 \triangleleft {V}_{1}$, where ${V}_{1}$ is seen as an abelian group and the action is given by $v_0\cdot v_1:=v_0+\dM v_1,$ for all $v_0\in V_0, v_1\in V_1$.

\begin{defi} Let $\mathbb V$ be a $2$-vector space. If $(\V,\Pi,\phi)$  is a quasi-Poisson groupoid such that the bivector field $\Pi$ and the trisection $\phi$ are both linear, we call $(\V,\Pi,\phi)$ a {\bf Lie-quasi-Poisson groupoid}.
\end{defi}

Consider the 2-vector space $\mathfrak{g}^*_{-1}\oplus{\g^*_{0}}\rightrightarrows{\mathfrak{g}^*_{-1}}$  given by the 2-term chain complex $
\mathfrak{g}^*_{0}\stackrel{l^*_1}\longrightarrow{\mathfrak{g}^*_{-1}}$, which is the dual of the complex $\g_{-1}\stackrel{l_1}{\longrightarrow}\g_0$. Consequently, by Theorem \ref{thm:int}, we generalize the fact that there is a Lie-Poisson structure on the dual space of a Lie algebra to the case of Lie 2-algebras.
\begin{thm}
The dual of a Lie $2$-algebra is a Lie-quasi-Poisson groupoid.
\end{thm}

We obtain our result through the integration of the Lie quasi-bialgebroid associated to a Lie 2-algebra. One can also obtain a quasi-Poisson Lie 2-group from a Lie 2-algebra via the integration of Lie 2-bialgebras. See \cite{chenxu} for more details.

\section*{Acknowledgement} We give our warmest
thanks to Jim Stasheff and the referee for very helpful comments that improve the paper. We also give our special thanks to Noriaki Ikeda, Zhangju Liu, Pavol Ševera and Chenchang Zhu for very
useful comments and discussions. We thank Marco Zambon for pointing out the reference \cite{Mehta} to us.

\footnotesize{

}


\begin{thebibliography}{999}
%\cite{Alexandrov:1995kv}

\bibitem{Anton1}
Alekseev, A. and Kosmann-Schwarzbach, Y.: Manin pairs and moment maps. \emph{J. Diff. Geom.} {\bf 56} (2000), 133-165.
%\yh{at the end of the title, we use ``.''. Furthermore, please xiaohong change the style of the author. Family name first, this is according to the style of this jourle. Follow the first, which I have revised.}

\bibitem{Anton2}
Alekseev, A., Kosmann-Schwarzbach, Y. and Meinrenken, E.: Quasi-Poisson manifolds. \emph{Canadian J. Math.}  {\bf 54} (2002), 3-29.

\bibitem{AKSZ}
Alexandrov, M., Kontsevich,  M., Schwartz, A.  and Zaboronsky, O.:
The geometry of the master equation and topological quantum field theory.
\emph{Int. J. Mod. Phys.}  A {\bf 12} (1997), 1405.


\bibitem{baez:2algebras}
 Baez,  J. C. and Crans, A. S.: Higher-Dimensional Algebra VI: Lie
 2-Algebras. \emph{Theory and Appl.  Categ.}  {\bf12} (2004),
 492-528.



\bibitem{baez:classicalstring}
Baez,  J. C., Hoffnung, A. E. and  Rogers, C. L.:  Categorified symplectic geometry and the classical string.
 \emph{Comm. Math. Phys.}, {\bf 293} (2010), no. 3, 701-725.
\bibitem{baez:string}
Baez,  J. C. and  Rogers, C. L.:
Categorified symplectic geometry and the string {L}ie 2-algebra.
\emph{Homology, Homotopy Appl.}, {\bf12 }(2010), no. 1,  221-236.

 \bibitem{Poncin} Bonavolont${\rm\grave{a}}$, G. and  Poncin, N.:
On the category of Lie $n$-algebroids. \emph{J. Geom. Phys.} {\bf 73} (2013), 79-90.

\bibitem{Bruce}
Bruce, A.: From $L_\infty$-algebroids to higher Schouten/Poisson structures. \emph{Rep. Math. Phys.} {\bf 67} (2011), no. 2, 157-177.

\bibitem{ADR} Cattaneo, A. S., Fiorenza, D. and Longoni, R.: Graded Poisson Algebras. \emph{Encyclopedia of Math. Phys.}, {\bf 2} (2006),  560-567.

\bibitem{AG} Cattaneo,  A. S. and  Felder, G.:  Relative formality theorem and quantisation of coisotropic submanifolds. \emph{Adv. Math.}, {\bf 208} (2007), no. 2, 521-548.

\bibitem{CZ}
Cattaneo, A. S. and Zambon, M.: A supergeometric approach to Poisson reduction. \emph{Comm. Math. Phys.} {\bf 318} (2013), no. 3, 675-716.

\bibitem{chenxu}
Chen, Z., Sti\'enon, M. and Xu, P.: Poisson 2-groups. \emph{J. Diff. Geom.} {\bf 94} (2013), no. 2, 209-240.

%\bibitem{street}
%Day, B. and Street, R.: Monoidal bicategories and Hopf algebroids. \emph{Adv. Math.}  {\bf129} (1997), 99-157.

\bibitem{Marco}
Fregier, Y. and  Zambon, M.: Simultaneous deformations and Poisson geometry.  \emph{Compos. Math.} 151 (2015), no. 9, 1763-1790.

%\bibitem{Grut}
%M. Grutzmann. $H$-twisted Lie algebroids, \emph{J. Geom. Phys.}, 2011, 61(2): 476-484.

\bibitem{strobl}
Hansen, M. and Strobl, T.: First Class Constrained Systems and Twisting of Courant
Algebroids by a Closed $4$-form. \emph{Fundamental interactions, {\rm115–144}, World Sci. Publ., Hackensack, NJ,} 2010.

\bibitem{Xu} Iglesias Ponte, D., Laurent-Gengoux, C. and Xu, P.:  Universal lifting theorem and quasi-Poisson groupoids.
 \emph{J. Eur. Math. Soc. (JEMS)} {\bf14} (2012), no. 3, 681-731.

%\bibitem{Ikeda}
%N. Ikeda and K. Uchino, QP-structures of degree 3 and 4D topological field theory. \emph{Comm. Math. Phys.} 303 (2011), no. 2, 317–330.

\bibitem{IkedaXu}
Ikeda, N. and Xu, X.: Canonical functions, differential graded symplectic pairs in supergeometry, and Alexandrov-Kontsevich-Schwartz-Zaboronsky sigma models with boundaries. \emph{J. Math. Phys},  {\bf 55} (2014), 113505.


\bibitem{VoronovHigherP}
Khudaverdian, H. M.  and Voronov, Th.: Higher Poisson brackets and differential forms. \emph{Geometric methods in physics,} 203-215, AIP Conf. Proc., 1079, \emph{Amer. Inst. Phys., Melville, NY,} 2008.

 \bibitem{Schwarzbach1}
Kosmann-Schwarzbach, Y.:  Jacobian quasi-bialgebras and quasi-Poisson Lie groups.
  \emph{Mathematical aspects of classical field theory (Seattle, WA, 1991), {\rm 459–489, Contemp. Math., 132,} Amer. Math. Soc., Providence, RI,} 1992.

 \bibitem{Schwarzbach2}
 Kosmann-Schwarzbach, Y.: Quasi, twisted, and all that… in Poisson geometry and Lie algebroid theory. \emph{The breadth of symplectic and Poisson geometry, {\rm 363-389,  Progr. Math., 232,} Birkh{\"a}user Boston, Boston, MA,} 2005.

  \bibitem{Schwarzbach3}
Kosmann-Schwarzbach, Y.: Poisson and symplectic functions in Lie algebroid theory. \emph{Higher structures in geometry and physics, {\rm 243-268, Progr. Math., 287,} Birkh{\"a}user/Springer, New York,} 2011.

  %\bibitem{Schwarzbach4}
%Y. Kosmann-Schwarzbach, Courant algebroids. A short history. \emph{SIGMA Symmetry Integrability Geom. Methods Appl.} 9 (2013), Paper 014, 8 pp.

 \bibitem{LadaMarkl}  Lada, T. and Markl, M.: Strongly homotopy Lie
algebras.
\emph{Comm. Algebra} {\bf 23} (1995),  no. 6, 2147-2161.

\bibitem{stasheff:introductionSHLA} Lada, T. and Stasheff, J.: Introduction to sh Lie algebras for
physicists. \emph{Int. Jour. Theor. Phys.}  {\bf 32} (1993), no. 7, 1087-1103.

%\bibitem{LiuShengXu}Z. Liu, Y. Sheng and X. Xu, Pre-Courant Algebroids and Associated Lie 2-Algebras, arXiv:1205.5898.


\bibitem{LWXmani}
Liu,  Z.,  Weinstein, A.  and Xu, P.: Manin triples for Lie
bialgebroids. \emph{J. Diff. Geom.}  {\bf 45} (1997), 547-574.
\bibitem{Mackenzie}
Mackenzie,  K. C. H.: Ehresmann doubles and Drinfel'd doubles for Lie algebroids and Lie bialgebroids.  \emph{J. Reine Angew. Math.}  {\bf 658} (2011), 193-245.
\bibitem{MackenzieXu}
 Mackenzie, K. C. H. and Xu, P.: Lie bialgebroids and Poisson groupoids. \emph{Duke Math. J.}  {\bf 73} (1994), no. 2, 415-452.


\bibitem{Mehta}
Mehta,  R. A.: On homotopy Poisson actions and reduction of symplectic $Q$-manifolds. \emph{Diff. Geom. Appl.}  {\bf 29} (2011), no. 3, 319-328.

\bibitem{Roytenberg1}
Roytenberg, D. and Weinstein, A.: Courant algebroids and strongly homotopy Lie
algebras. \emph{Lett. Math. Phys.},  {\bf 46} (1998),  no. 1, 81-93.

\bibitem{Roytenberg2} Roytenberg,  D.: On the structure of graded symplectic supermanifolds and Courant algebroids. \emph{Quantization, Poisson Brackets and Beyond}, 169-185,  Contemp. Math., 315, \emph{ Amer. Math. Soc., Providence, RI,} 2002.

    \bibitem{Roytenberg3}
Roytenberg, D.:
Quasi-Lie bialgebroids and twisted Poisson manifold.
\emph{Lett. Math. Phys.} {\bf 61} (2002), no. 2, 123-137.

\bibitem{Roytenberg4}
Roytenberg, D.: \emph{Courant algebroids, derived brackets and even
symplectic supermanifolds}. Ph. D thesis, UC Berkeley, 1999,
arXiv:math.DG/9910078.

\bibitem{Roytenberg5}Roytenberg, D.: On weak Lie $2$-algebras. \emph{ XXVI Workshop on Geometrical Methods in Physics.} 180-198.   AIP Conf. Proc., 956. \emph{Amer. Inst. Phys., Melville, NY,} 2007.

\bibitem{Schatz}
Sch${\rm\ddot{a}}$tz,  F.: \emph{Coisotropic submanifolds and the BFV-complex}. Ph. D thesis, University Z${\rm\ddot{u}}$rich, 2009


\bibitem{Severa}
Ševera,  P.: Some title containing the words ``homotopy'' and ``symplectic'', e.g. this one. \emph{Travaux math$\acute{e}$matiques. Fasc. XVI,} 121-137, Trav. Math., XVI, \emph{ Univ. Luxemb., Luxembourg,} 2005.
\bibitem{Severa1}
Ševera, P. and Weinstein,  A.: Poisson geometry with a 3-form background. \emph{Prog. Theor. Phys.
Suppl.}, {\bf 144} (2001), 145-154.



%\bibitem{LiuSheng}Y. Sheng and Z. Liu, Lebniz 2-alegbras and twisted Courant algebroids. \emph{Comm. Algebra.} 41 (05),1929-1953.

 \bibitem{shengzhu} Sheng, Y. and  Zhu, C.:  Higher Extensions of Lie Algebroids.  \emph{Comm. Contemp. Math.} (2016), DOI: 10.1142/S0219199716500346.

 \bibitem{shengzhu1} Sheng, Y. and  Zhu, C.: Semidirect products of representations up to homotopy. \emph{Pacific J. Math.}, {\bf 249} (2011),  no. 1, 211-236.

 \bibitem{shengzhu2} Sheng, Y. and  Zhu, C.:
Integration of semidirect product Lie 2-algebras. \emph{Int. J. Geom. Methods Mod. Phys.}  {\bf 9}  (2012),  no. 5, 1250043.

\bibitem{stasheff:shla}
Stasheff, J.:
 Differential graded {L}ie algebras, quasi-Hopf algebras and higher
  homotopy algebras.   \emph{Quantum groups (Leningrad, 1990)}, 120-137,
  Lecture Notes in Math., 1510, \emph{ Springer, Berlin,} 1992.
 \bibitem{Voronov0}
Voronov, Th.: Graded manifolds and Drinfeld doubles for Lie bialgebroids. \emph{Contemp. Math.} {\bf 315} (2002), 131-168.
\bibitem{Voronov1}
Voronov,  Th.: Higher derived brackets and homotopy algebras. \emph{J. Pure Appl. Algebra}  {\bf 202} (2005), no. 1-3, 133-153.

\bibitem{Alan}
Weinstein, A.: Omni-Lie algebras. Microlocal analysis of the Schrödinger equation and related topics (Japanese) (Kyoto, 1999). \emph{S${\bar{u}}$rikaisekikenky${\bar{u}}$sho K${\bar{u}}$ky${\bar{u}}$roku} {\bf 1176} (2000), 95-102.

\bibitem{XuM}
Xu, X.: Twisted Courant algebroids and coisotropic Cartan geometries. \emph{J. Geom. Phys.}  {\bf 82} (2014), 124-131.


\end{thebibliography}
 \end{document}